\documentclass[12pt]{amsart}
\usepackage{amssymb}
\usepackage{amsmath}
\usepackage{amsthm}
\usepackage{color}
\usepackage{amsfonts}
\usepackage{wrapfig}
\usepackage{graphicx}
\usepackage{float}
\usepackage{caption}
\usepackage{subcaption}
\usepackage[rightcaption]{sidecap}
\usepackage[section]{placeins}
\usepackage[font=small,labelfont=bf]{caption}

\usepackage{url}

\newtheorem{theorem}{Theorem}[section]
\newtheorem{lemma}[theorem]{Lemma}

\newtheorem{definition}{Definition}

\newtheorem{example}{Example}

\theoremstyle{remark}
\newtheorem{remark}[theorem]{Remark}
\newtheorem*{notation}{Notation}

\numberwithin{equation}{section}

\DeclareMathOperator{\supp}{supp}

\newcommand{\R}{\mathbb{R}}

\newcommand{\Fo}{\mathcal{F}_c^0(\R^+)}
\newcommand{\Ft}{\mathcal{F}_c(\R^+)}
\newcommand{\Fth}{\mathcal{F}_c(\R)}

\title{Arithmetic of fuzzy numbers and intervals - \\
\tiny  a new perspective with examples}
\author[Jan Schneider]{Jan Schneider}
\thanks{\textit{Special thanks for general endorsement:} Sweetchris}

\address{Faculty of Computer and Management Science,
Wroclaw University of Technology,
ul. Smoluchowskiego 25,
50-371 Wroclaw, Poland}
\email{jan.schneider@pwr.wroc.pl}
\subjclass[2010]{}
\begin{document}

\begin{abstract}
This article is meant to give a lucid and widely accessible, self-contained account of a novel way of performing arithmetic operations on fuzzy intervals. Based on two formulae of generalized inversion (the first in close analogy to the inversion of cumulative distribution functions in probability and statistics, and subsequent re-inversion which in this form seems to be new in the literature) the technique could prove to be of some importance for the practical handling of fuzzy intervals whose characterizing functions are discontinuous and/or not strictly monotonic, e.g.~piecewise constant.
Notwithstanding its innovative nature, the article is written in the style of an introductory textbook featuring a multitude of illustrated examples.
\end{abstract}

\maketitle

\setcounter{section}{-1}

\section{Summary and Structure}
This article guides the reader through four consecutive definitions of widening scope. There are thus four sections: \smallskip\par\noindent
\textbf{Section~1:} A first simplified definition of fuzzy numbers, as just members of a set without (arithmetic) structure, is given. The characterizing functions are continuous and their support is assumed to be a closed interval within the positive reals $\R^+$.
\par
The arithmetic operations ``$\oplus, \odot$'' of addition and multiplication are defined. The formulation is mathematically equivalent to Lotfi A. Zadeh’s original definition \cite{Zadeh65} of binary operations via two-dimensional fuzzy vectors and the extension principle, but closer to the more user-friendly definitions of R. Goetschel and W. Voxman~\cite{GV83, GV86}, albeit introducing different terminology and notation.\par\noindent
\textbf{Section~2:} Starting with section 2 we prefer to speak more generally about fuzzy intervals regarding fuzzy numbers a subcase.
To make possible the progression from \textit{Definition~1} (continuous characterizing functions and support on the positive real line) to \textit{Definition~2} (semi-continuous characterizing functions, support on the positive reals) the notion of generalized inverse functions as known in  probability and statistics is briefly discussed and illustrated. An explicit formula for the converse operation is given. The consequent  definitions of fuzzy arithmetic are the main theoretical contribution of this article. \par\noindent
The representation of a fuzzy interval as an ordered pair of functions is adapted.\par\noindent
\textbf{Section~3} extends the definitions of sections~1 and~2 to the case of characterizing functions of mixed support (that is containing~0).\par\noindent
\textbf{Section~4:} In the concluding section the parametric representation of fuzzy intervals as defined in \cite{GV83, GV86} and arithmetic operations on them is laid out, and it is argued that the parametric representation be used as the principal, primary and basic one, as has been the trend in recent publications.
\tableofcontents
\newpage

\section{$\Fo$}
\markboth{Arithmetic of Fuzzy Intervals and Numbers}{Jan Schneider}

\begin{definition}\label{def1}
For the purpose this section:
A fuzzy number $\xi$ is defined by, and identified with a continuous function (called the {\it characterizing function} of that fuzzy
number) $\xi(\cdot):\R^+\to[0,1]$ which has the following properties:
\begin{itemize}
\item[(1)] compact, positive support: the closure of the set of points $x\in\R^+$, where $\xi(x)>0$ is a closed interval $[l,r]\subset\R^+,$ for some $l,r\in\R^+$
\item[(2)] there is exactly one point $m\in(l,r)$ such that $\xi(m)=1,$
\item[(3)] $\xi(\cdot)$ is strictly increasing on the interval $[l,m],$
\item[(4)] $\xi(\cdot)$ is strictly decreasing on the interval  $[m,r].$
\end{itemize}
\end{definition}
\noindent Notice that condition~(1) of \textit{Definition~1} implies that $\xi_l(l)=0$ and $\xi_r(r)=0.$\par
\smallskip
\noindent The collection of fuzzy numbers $\xi$ satisfying \textit{Definition~\ref{def1}} shall be denoted by $\Fo$.\par
\smallskip
\noindent\textit{Convention:} Throughout the article when writing $\R^+$ we mean that $0\in\R^+$.\par
\smallskip
\noindent For future computational comfort, the increasing (resp. the decreasing) components of $\xi$ are singled out and treated as separate
entities denoted by $\xi_l$ (resp. $\xi_r$). (As in left, right.)\par\smallskip
For convenience we may now rewrite $\xi$:\par
\smallskip
\noindent\textbf{Definition 1'.}
\begin{equation}\label{lrform1}
\xi(x)=
\begin{cases}
\xi_l(x)&\text{ for $x\in[l,m],$}\\
1 &\text{ for $x=m,$}
\\
\xi_r(x)&\text{ for $x\in[m,r].$}
\end{cases}
\end{equation}
with $\xi,\xi_l,\xi_r$ meeting conditions (1)-(4) of \textit{Definition~1}.\par\noindent\vspace{0.1mm}
\noindent We shall also sometimes write $\xi_m$ (as in middle) for the one point $x_o = m$, where $\xi(m)=1$.\par
\noindent
(The formally faulty overkilling of $m$ in the definition of $\xi$ above, using closed intervals instead of half-closed is intentional and will find its justification at a later point in this article. We need to have $\xi_l$ and $\xi_r$ defined on closed intervals, as a formal requirement for the process of inverting back and forth)\par
\smallskip\noindent
\textbf{Definition 1''}
Equivalently one may think of a fuzzy number as the ordered pair $\bigl(\xi_l(\cdot),\xi_r(\cdot)\bigr)$, with $m$ again being implicit from \mbox{$\xi_l(m)=\xi_r(m)=1$.}\par

\bigskip
\noindent\textbf{Terminology:}\par\noindent
For a given fuzzy number $\xi$ the functions \mbox{$\xi_l:[l,m]\to[0,1]$} and \mbox{$\xi_r:[m,r]\to[0,1]$} shall be for the purpose of this article called the left and right \textit{fuzzy endpoints} of $\xi$.\par\noindent

\smallskip
We may colloquially say:\textit{``Here we have a real number $``m"$ which is fuzzy to the left with $\xi_l(x)=....$ and fuzzy to the right with
$\xi_r(x)=....$"}\par\noindent
\smallskip

\begin{figure}[H]
\begin{center}
\includegraphics[height=12mm]{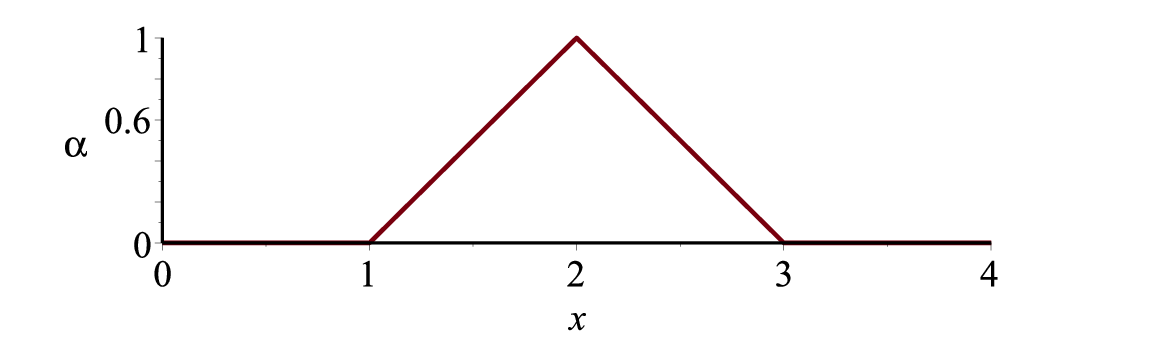}
\includegraphics[height=14mm]{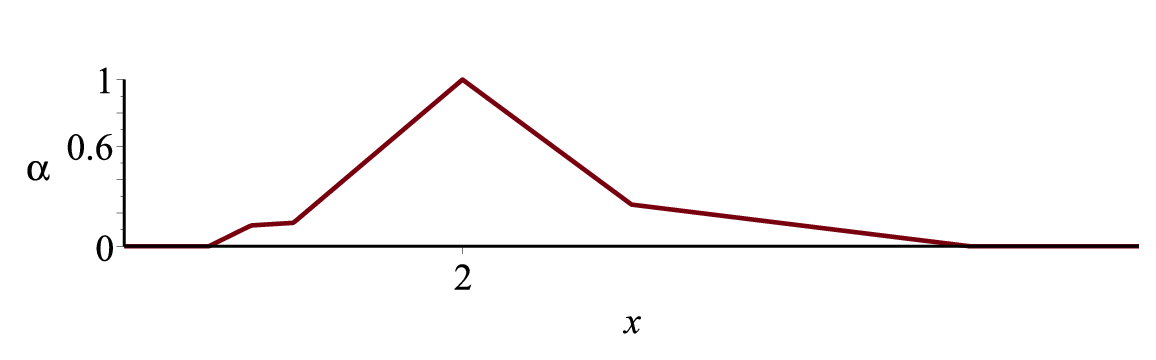}
\includegraphics[height=12mm]{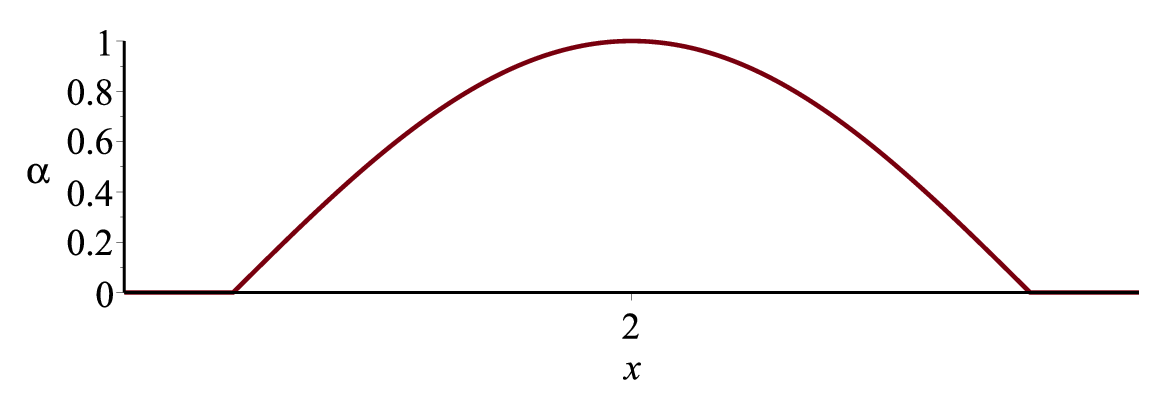}
\caption{Examples of fuzzy numbers: (a) $\xi_l$, $\xi_r$ -- linear (b) $\xi_l$, $\xi_r$ -- piecewise linear, (c)
$\xi$ -- of sine type }\label{plr}
\end{center}
\end{figure}

\noindent Below, in Fig.~\ref{plr} several examples of "fuzzy 2-s" as per \textit{Definition~1} are charted:\par
\begin{figure}[htbp]
\begin{center}
\includegraphics[width=40mm]{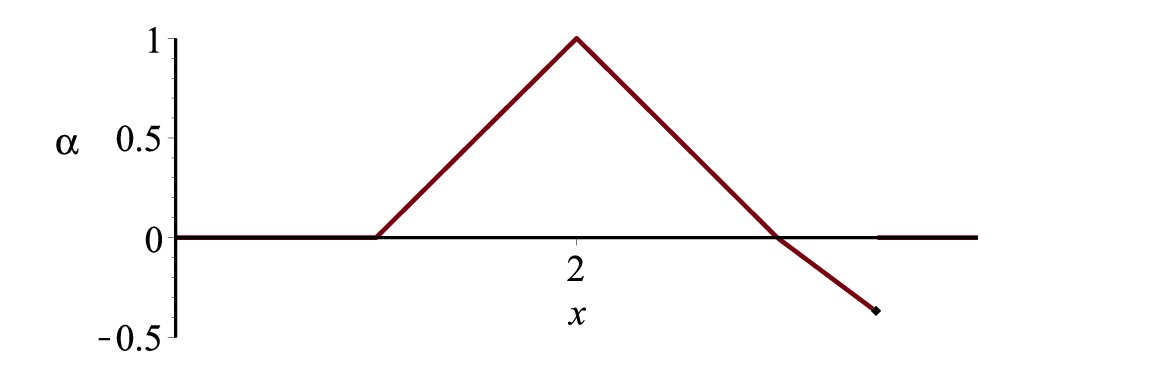}
\includegraphics[width=55mm]{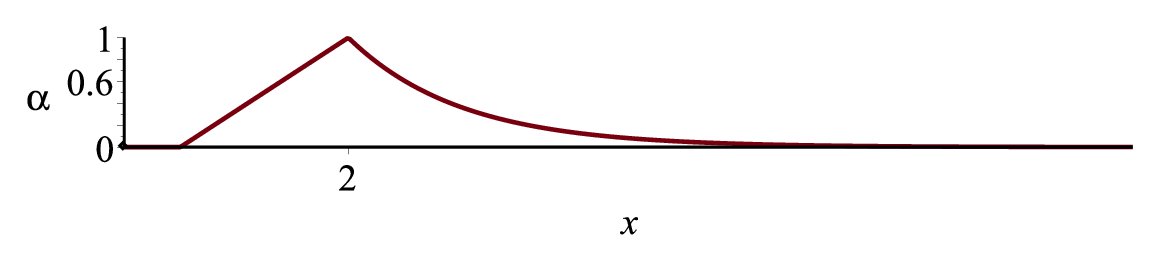}
\includegraphics[width=40mm]{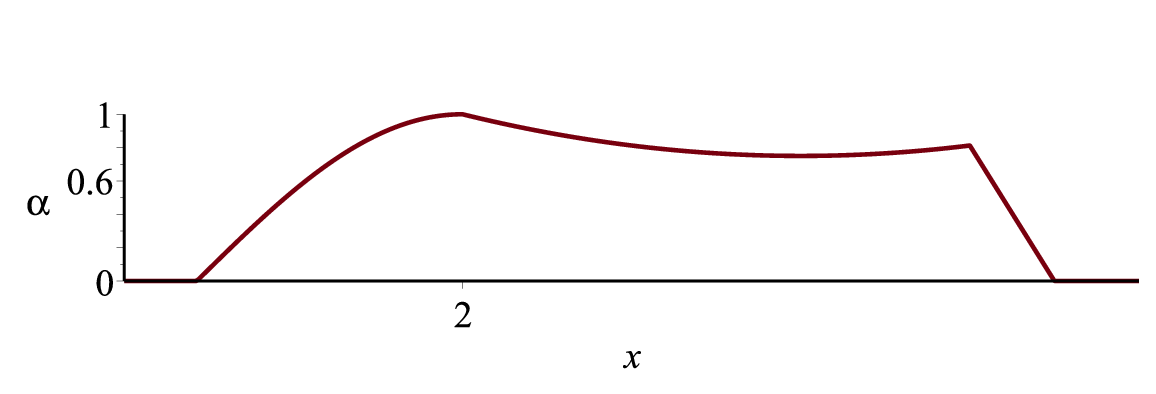}
\caption{Functions which do not characterize fuzzy numbers: (a) function with negative values, (b) function with unbounded
support, (c) monotonicity condition not fulfilled.}\label{nlr}
\end{center}
\end{figure}

\noindent On the other hand: the functions illustrated in Fig.~\ref{nlr} do not characterize fuzzy numbers in the sense of \textit{Definition~\ref{def1}}.
In Fig.~\ref{nlr} (a) the function is not non-negative. In Fig.~\ref{nlr} (b) the set of points, where the
function is positive is not contained in an interval of finite length. Finally, the function in Fig.~\ref{nlr} (c)
does not satisfy the monotonicity conditions (3), (4) of \textit{Definition~\ref{def1}}.

\begin{remark} In the scientific community dealing with fuzzy numbers it is customary to denote the values of characterizing
functions  by $\alpha$ (not e.g. $y$), where \mbox{$0\leq\alpha\leq1$}. When analyzing graphs of fuzzy numbers one needs to be ever conscientious
about whether and to what degree the ordinate axis is scaled.
(We shall later see that $\supp[\xi\odot\eta]=\supp[\xi]\cdot\supp[\eta]$.) In this presentation we shall look at ``small'' numbers and as a rule employ non-scaled graphs.
\end{remark}

\subsection{Addition and multiplication on $\Fo$}
\begin{remark}
Because of the monotonicity conditions (1.3), (1.4) of \textit{Definition~1} the fuzzy components $\xi_l$ and $\xi_r$ of $\xi$ are invertible with  \mbox{$\xi_l^{-1}:[0,1]\to\R$} with \mbox{$\xi_l^{-1}(\R)=[l,m]$} and \mbox{$\xi_r^{-1}:[0,1]\to\R$} with \mbox{$\xi_r^{-1}(\R)=[m,r]$}, which are again invertible and those new inverses are once more strictly monotone and have values $\alpha\in[0,1]$ and we may write:
\end{remark}
\begin{subsubsection}
{Addition} of fuzzy numbers is defined in the following way:
\begin{equation}\label{dodawanie}
\begin{split}
(\xi\oplus\eta)_l(x)\,= \, (\xi_l^{-1}+\eta_l^{-1})^{-1}(x),\\
(\xi\oplus\eta)_r(x)\,= \,(\xi_r^{-1}+\eta_r^{-1})^{-1}(x).
\end{split}
\end{equation}

\smallskip
This implies that $(\xi\oplus\eta)(m)=1$, where $m=\:\xi^{-1}(1)+\eta^{-1}(1)$, so using the notation of \textit{Definition~1'} we have
\begin{equation}
\mbox{$(\xi\oplus\eta)_m=\xi_m+\eta_m$}.
\end{equation}
\par

\smallskip
\noindent In interval notation (\textit{Definition~1''}) we write:
\begin{equation}
(\xi_l,\xi_r)\oplus(\eta_l,\eta_r)=\biggl((\xi_l^{-1}+\eta_l^{-1})^{-1}\,\textbf{,}\,(\xi_r^{-1}+\eta_r^{-1})^{-1}\biggr)
\end{equation}

\begin{remark} The operation given in (\ref{dodawanie}) is well defined since being the sum of strictly monotone functions $(\xi_l^{-1}+\eta_l^{-1})$ and
$(\xi_r^{-1}+\eta_r^{-1})$ are also strictly monotone and therefore are invertible on their natural domain $[0,1]$.
\end{remark}
\end{subsubsection}
\smallskip
\begin{subsubsection}{Multiplication} is defined analogically
\begin{equation}\label{mnozenie}
\begin{split}
(\xi\odot\eta)_l(x)\,= \,(\xi_l^{-1}\cdot\eta_l^{-1})^{-1}(x),\\
(\xi\odot\eta)_r(x)\,= \,(\xi_r^{-1}\cdot\eta_r^{-1})^{-1}(x).
\end{split}
\end{equation}
With
\begin{equation}
(\xi\odot\eta)(m)\,=\: 1\, \quad\text{where \:$m=\:\xi^{-1}(1)\cdot\eta^{-1}(1)$},
\end{equation}
or
\begin{equation}
(\xi\odot\eta)_m=\xi_m\cdot\eta_m.
\end{equation}
\smallskip
and again
\begin{equation}
(\xi_l,\xi_r)\odot(\eta_l,\eta_r)=\left((\xi_l^{-1}\cdot\eta_l^{-1})^{-1},(\xi_r^{-1}\cdot\eta_r^{-1})^{-1}\right)
\end{equation}
\end{subsubsection}

The formulae for addition and multiplication of fuzzy numbers given above may be explicated in a convenient albeit quite
informal way:\par
\smallskip
\noindent
\textit{A fuzzy $5$ times a fuzzy $3$ is a fuzzy $15$ and their sum is a fuzzy $8$. The tricky part is to compute the left and right fuzzy
parts:
Here is what you do:}\par\smallskip
{\it Draw the graphs of the two fuzzy numbers on a piece of paper acting as the the coordinate plane of $x$ and $\alpha$. Rotate the piece of paper by
90~degrees to interchange the roles of $x$ and $\alpha$. Add or multiply in standard manner as functions the
components of the rotated functions. Rotate back the coordinate plane to see the outcome.}

\begin{subsection}{Algebraic properties of $\Fo$}
It follows directly from the definition that the operations of addition ``$\oplus$'' and multiplication
``$\odot$'' are commutative and associative, i.e.
for arbitrary $\xi,\eta,\mu\in\Fo,$
\begin{equation*}
\begin{split}
\xi\oplus\eta&=\eta\oplus\xi\quad\text{ and }\quad\xi\odot\eta=\eta\odot\xi.\\
(\xi\oplus\eta)\oplus\mu&=\eta\oplus(\xi\oplus\mu)\quad\text{ and }\quad(\xi\odot\eta)\odot\mu=\eta\odot(\xi\odot\mu).
\end{split}
\end{equation*}
Moreover, the standard distributive property
\begin{equation}\label{distributivity}
\xi\odot(\eta\oplus\mu)=(\xi\odot\eta)\oplus(\xi\odot\mu)
\end{equation}
is satisfied as well. (for (\ref{distributivity}) to hold, it is essential that the support of the characterizing functions be contained in the positive reals. Distributivity will not hold for
the extension of section~3).
\end{subsection}

We shall close this section with several examples illustrating the operations of addition and multiplication.
\begin{example}
The simplest examples of fuzzy numbers, and those that appear in the literature most often, are numbers whose
fuzzy endpoints are linear. (Such numbers are commonly referred to as ``triangle numbers" and in short denoted by $tr_{(l,m,r)}(\cdot)$.
So let us consider two such fuzzy numbers $\xi$ and $\eta$ defined by their characterizing functions:
\begin{equation}\label{xiexample1}
\xi(x)=\begin{cases}\xi_l(x)=x-1&\text{ for $x\in[1,2],$}\\
1&\text{ for $x=2$},\\
\xi_r(x)=3-x&\text{ for $x\in[2,3],$}
\end{cases}
\end{equation}
and
\begin{equation}\label{etaexample1}
\eta(x)=\begin{cases}\eta_l(x)=\frac{1}{2}x-\frac{5}{2}&\text{ for $x\in[5,7],$}\\
1&\text{ for $x=7$},\\
\eta_r(x)=-\frac{1}{3}x+\frac{10}{3}&\text{ for $x\in[7,10].$}
\end{cases}
\end{equation}

\begin{figure}[htbp]
\begin{center}
\includegraphics[width=120mm]{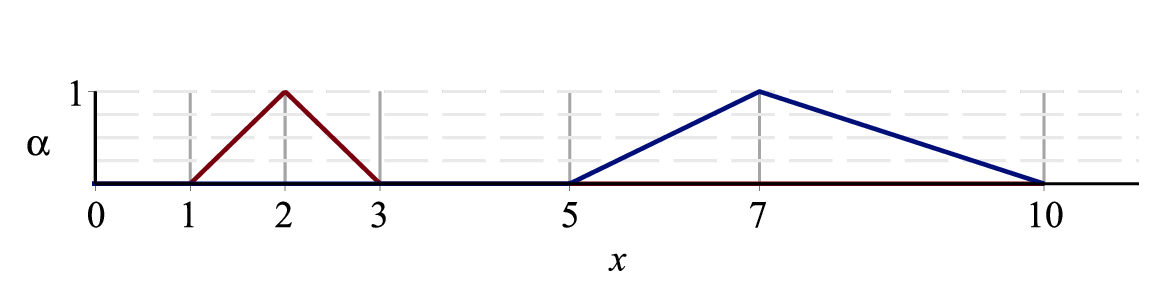}
\caption{Two fuzzy numbers $\xi$ and $\eta$ with linear fuzzy endpoints.}\label{xieta}
\end{center}
\end{figure}

Now whether we want to add or multiply $\xi$ and $\eta$ we need to first find the inverse functions $\xi_l^{-1},\xi_r^{-1}$ and
$\eta_l^{-1},\eta_r^{-1}\,$:\smallskip
\begin{align*}
\xi_l^{-1}(\alpha)=\,\,&\alpha+1,&\eta_l^{-1}(\alpha)=\,&2\alpha+5,\\
\xi_r^{-1}(\alpha)=\,\,&3-\alpha,&\eta_r^{-1}(\alpha)=\,&10-3\alpha.
\end{align*}
\noindent To perform the addition we use \eqref{dodawanie}. For $\alpha\in[0,1]$ we obtain
\begin{align*}
\xi_l^{-1}(\alpha)+\eta_l^{-1}(\alpha)=&\,3\alpha+6,\\
\xi_r^{-1}(\alpha)+\eta_r^{-1}(\alpha)=&\,13-4\alpha.
\end{align*}
Thus, inverting once again by \eqref{dodawanie} we arrive at
\begin{equation*}
\xi\oplus\eta\,(x)=
\begin{cases}
(\xi\oplus\eta)_l(x)=\frac{1}{3}x-2 &\text{ for $x\in[6,9],$}\\
1 &\text{ for $x=9,$}\\
(\xi\oplus\eta)_r(x)=-\frac14x+\frac{13}{4} &\text{ for $x\in[9,13].$}
\end{cases}
\end{equation*}
with $(\xi\oplus\eta)_m=\xi_m+\eta_m=2+7=9$.\par
\bigskip
\begin{figure}[htbp]
\begin{center}
\includegraphics[width=120mm]{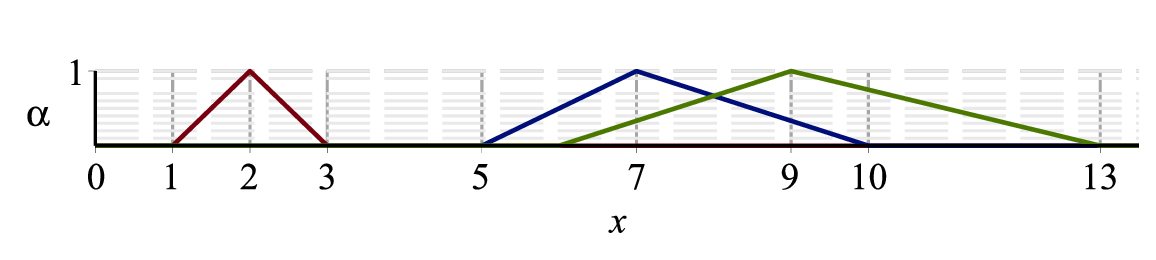}
\caption{Two fuzzy numbers $\xi$, $\eta$ and their sum $\xi\oplus\eta$}\label{xieta+}
\end{center}
\end{figure}
Fig.~\ref{xieta+} shows $\xi$, $\eta,$ and $\xi\oplus\eta.$ We note that the fuzzy addition operation~``$\oplus$''
preserves linearity of components.

\smallskip\noindent
Now let us move on to multiplication. Going by \eqref{mnozenie} we get
\begin{align*}
\xi_l^{-1}(\alpha)\cdot\eta_l^{-1}(\alpha)=\,\,&2\alpha^2+7\alpha+5,\\
\xi_r^{-1}(\alpha)\cdot\eta_r^{-1}(\alpha)=\,\,&3\alpha^2-19\alpha+30.
\end{align*}
for $\alpha\in[0,1].$
We now need to find the inverse functions of the obtained products. For the function
$\xi_l^{-1}(\alpha)\cdot\eta_l^{-1}(\alpha)$
we have $2\alpha^2+7\alpha+5=x.$ Since $\alpha$ changes between 0 and 1, we see that $x$ changes between $5$ and
$14.$ Thus we need to solve
the quadratic equation $2\alpha^2+7\alpha+5-x=0$ for $\alpha\in[0,1]$ and $x\in[5,14]$. We calculate the
discriminant $\Delta=9+8x$ and see that $\alpha=\frac{-7+\sqrt{9+8x}}{4}$ for
$x\in[5,14].$ The function $\xi_r^{-1}(\alpha)\cdot\eta_r^{-1}(\alpha)$ can be inverted in a similar way (the details
are left as an exercise for the reader). Clearly, $(\xi\odot\eta)_m=\xi_m\cdot\eta_m = 2\cdot7=14$, so finally we
obtain
\begin{equation*}
(\xi\odot\eta)(x)=
\begin{cases}
(\xi\odot\eta)_l(x)=\frac{-7+\sqrt{9+8x}}{4}&\text{ for $x\in[5,14],$}\\
1 &\text{ for $x=14,$}\\
(\xi\odot\eta)_r(x)=\frac{19-\sqrt{1+12x}}{6}&\text{ for $x\in[14,30].$}
\end{cases}
\end{equation*}

\begin{figure}[htbp]
\begin{center}
\includegraphics[width=125mm]{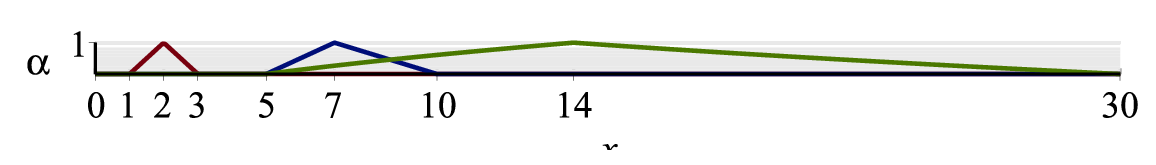}
\caption{Fuzzy numbers $\xi$, $\eta$ and their product $\xi\odot\eta$.}\label{xieta*}
\end{center}
\end{figure}

\begin{remark}
The product $\xi\odot\eta$ of two fuzzy numbers with linear characterizing functions has components that are
not linear, but are rather of a square root type.
\end{remark}
\end{example}

\begin{example}\label{sinus}
Let us consider the fuzzy number $\xi(x)=\chi_{[0,\pi]}(x)\cdot\sin x $ and decompose and write it in the form of \textit{Definition 1'}:
\begin{equation*}
\xi(x)=
\begin{cases}\sin x&\text{ for $x\in[0,\frac\pi2],$}\\
1&\text{ for $x=\frac\pi2$}\\
\sin x&\text{ for $x\in[\frac\pi2,\pi].$}
\end{cases}
\end{equation*}
We would like to calculate $\xi\odot\xi$, that is, $\sin^{\odot2}$.  We have $\xi_l(x)=\sin x$ for
$x\in[0,\frac{\pi}{2}]$ and $\xi_r(x)=\sin x$ on $[\frac\pi2,\pi].$ Thus $\xi_l^{-1}(\alpha)=\arcsin\alpha$ for
\linebreak $\alpha\in[0,1],$ and $\xi_r^{-1}(\alpha)=\pi-\arcsin(\alpha)$ for $\alpha\in[0,1].$
Both these functions are shown in Fig.~\ref{3sin}.
\begin{figure}[htbp]
\begin{center}
\includegraphics[width=60mm]{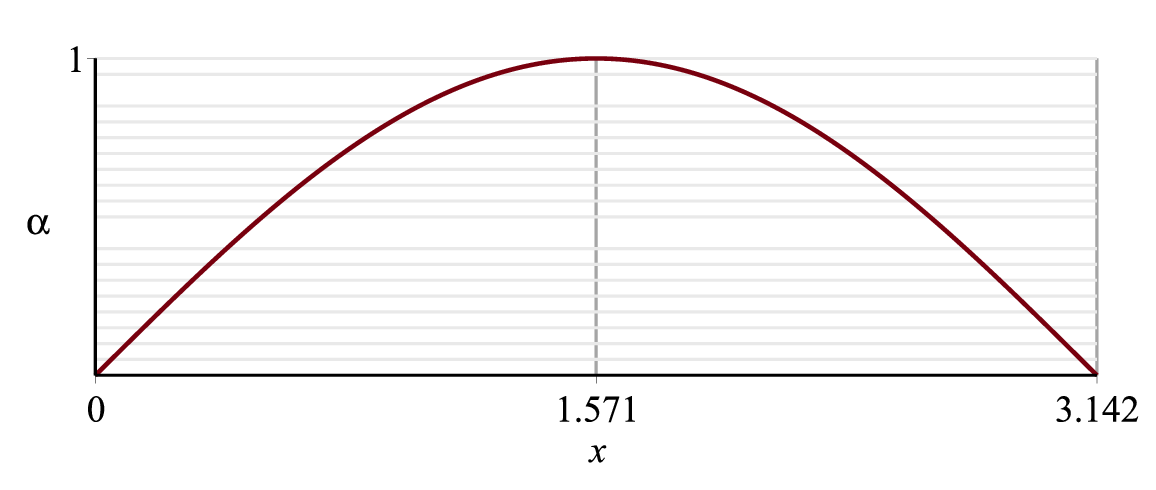}
\includegraphics[width=30mm]{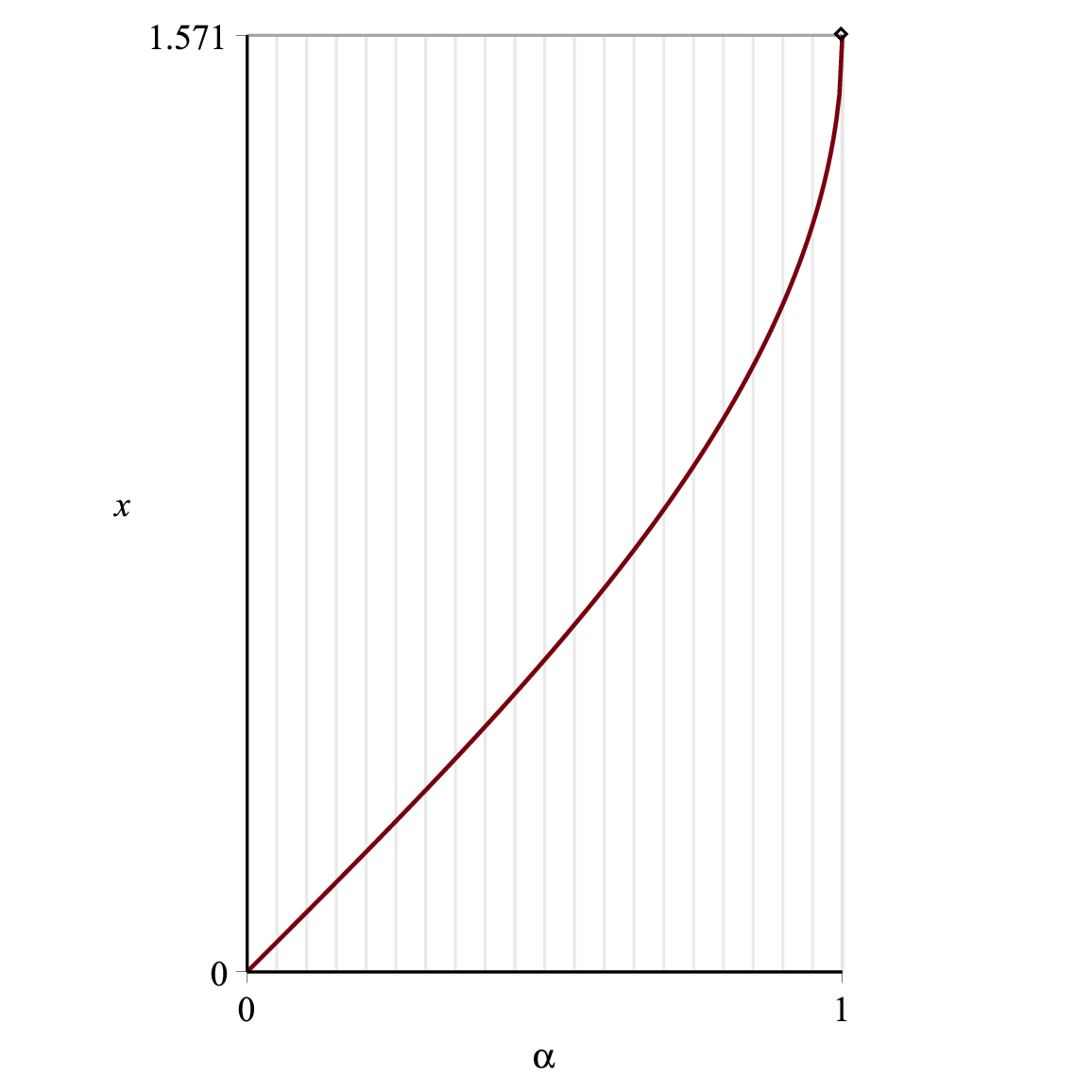}
\includegraphics[width=30mm]{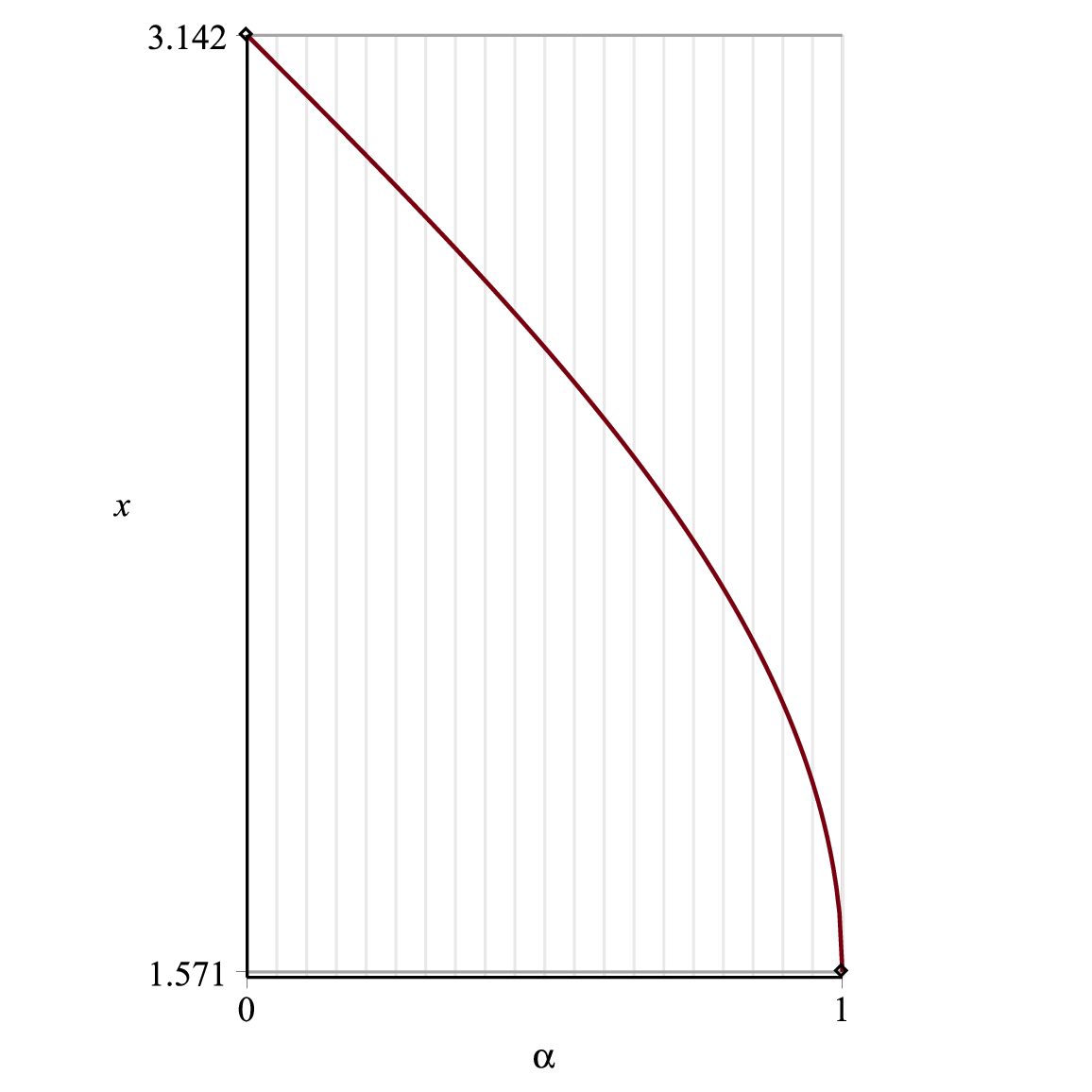}
\caption{(a) the function $\xi(x)=\chi_{[0,\pi]}(x)\sin x,$ (b) the inverse $\xi_l^{-1}(\alpha),$ (c) the inverse
$\xi_r^{-1}(\alpha).$}\label{3sin}
\end{center}
\end{figure}

\bigskip
\noindent In order to find $(\xi\odot\xi)_l$ we need to invert the function $x=\arcsin^2\alpha$. Similarly, for
$(\xi\odot\xi)_r$ we need to invert $x=(\pi-\arcsin(\alpha))^2.$ In this way we obtain
\begin{equation*}
\xi\odot\xi\,(x)=
\begin{cases}
(\xi\odot\xi)_l(x)=\sin(\sqrt{x})&\text{ for $x\in[0,\frac{\pi^2}{4}],$}\\
1 & \text { for $x=\frac{\pi^2}{4},$}\\
(\xi\odot\xi)_r(x)=\sin(\pi-\sqrt{x})=\sin(\sqrt x)&\text{ for $x\in[\frac{\pi^2}{4},\pi^2].$}
\end{cases}
\end{equation*}
\begin{figure}
\includegraphics[height=30mm]{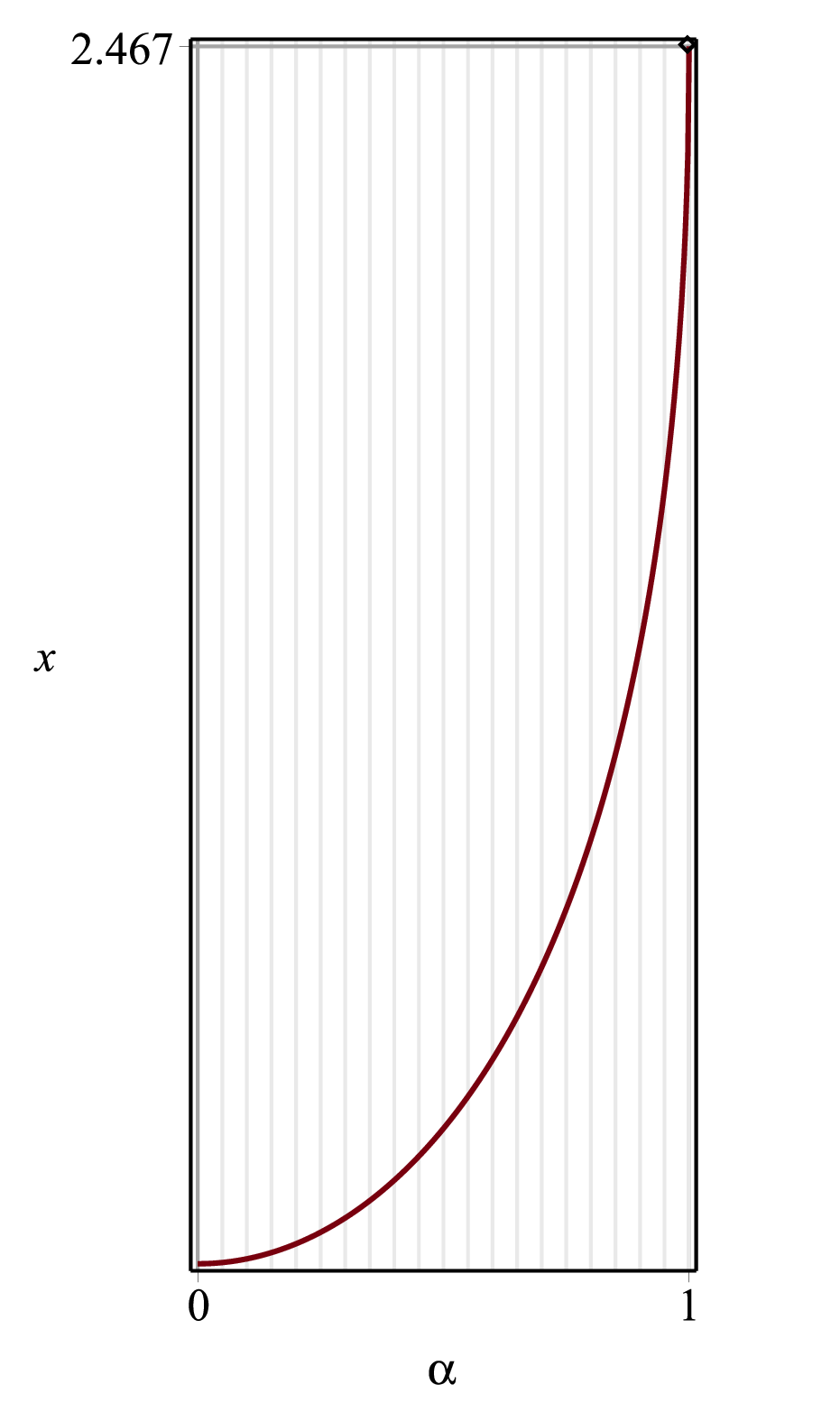}
\includegraphics[height=30mm]{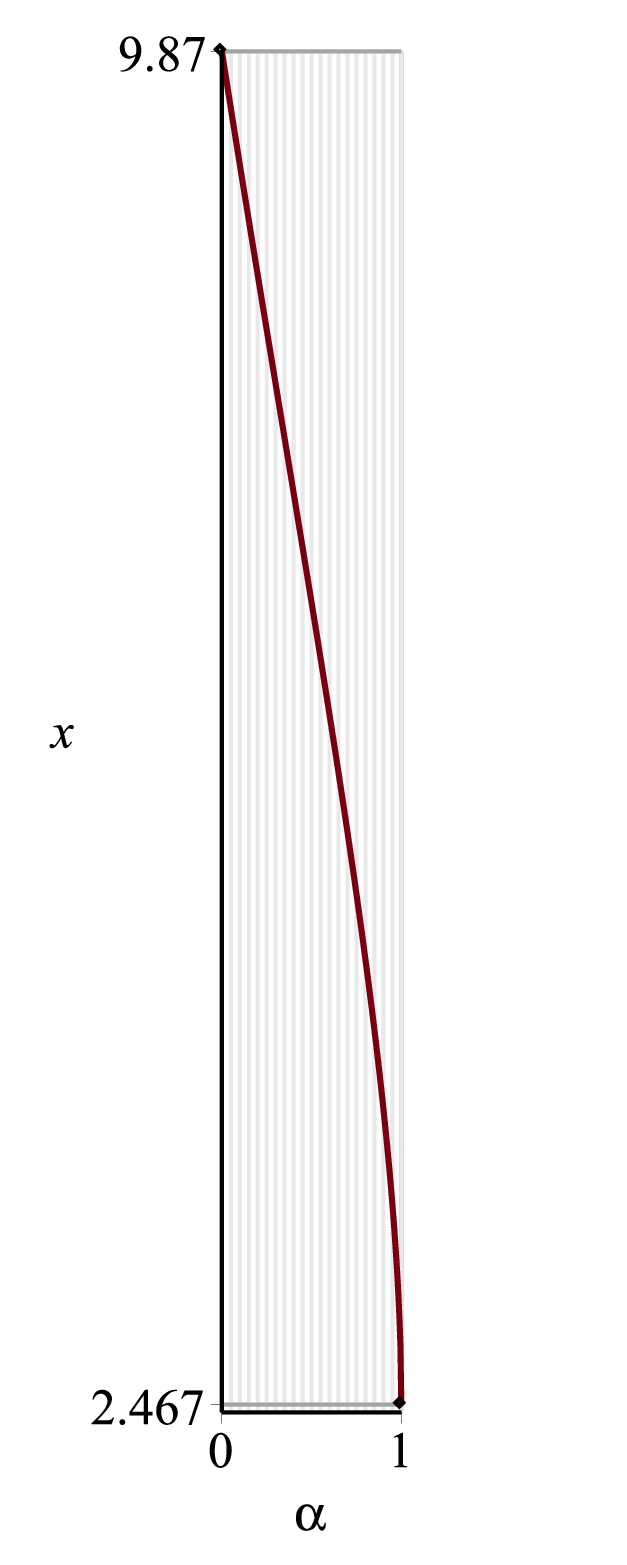}
\includegraphics[width=80mm]{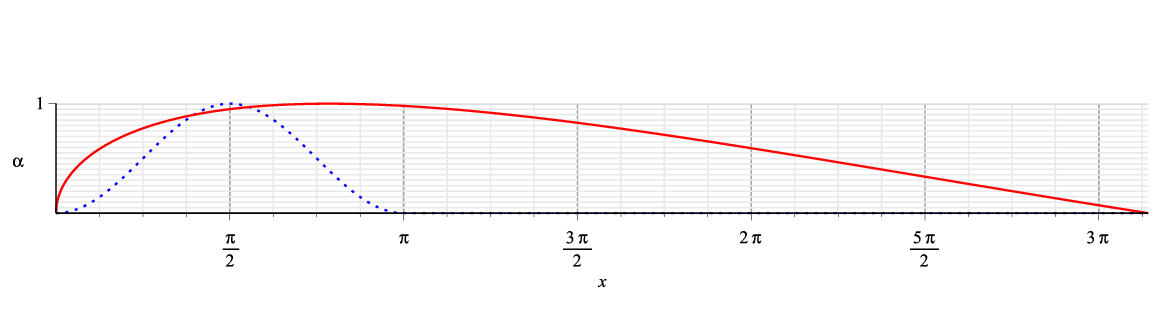}
\caption{$(\xi^{-1}_l)^2,(\xi^{-1}_r)^2$ and $\sin^{\odot2}$. For comparison (dotted blue line) $\sin^2$.}\label{sinsin}
\end{figure}
\begin{remark}
We see that in the realm (domain) of fuzzy numbers we have an identity $\sin^{\odot2}(x)=\sin(\sqrt x)$ with supports $[0,\pi]$ and $[0,\pi^2]$
respectively.\par\vspace{1mm}\noindent
- This is not a coincidence but an example indicating of how functions of fuzzy quantities work.
\end{remark}
\end{example}

\begin{example}
Let $\xi$ be as in the previous example, that is,  $\xi(x)=\chi_{[0,\pi]}(x)\sin x$ and $\eta(x)=\chi_{[\frac{3\pi}
2,\frac{5\pi} {2}
]}\cos x.$ For $\eta$
we have $\eta_l^{-1}(\alpha)=\arcsin(\alpha)+\frac{3}{2}\pi$ and $\eta_r^{-1}(\alpha)=\arcsin(-\alpha)+5\pi\slash
2$ for $\alpha\in[0,1].$ Thus,
we invert the functions $x=\arcsin(\alpha)(\arcsin\alpha+3\pi\slash 2)$ and
$x=(\arcsin(-\alpha)+\pi)(\arcsin(-\alpha)+5\pi\slash 2)$ to arrive at (see Fig. \ref{sincos})
\begin{equation*}
\xi\odot\eta =
\begin{cases}
(\xi\odot\eta)_l=\sin\left(-\frac{3}{4}\pi+\frac{1}{2}\sqrt{\frac{9}{4}\pi^2+4x}\right)&\text{ for
$x\in[0,\pi^2],$}\\
1 & \text{ for $x=\pi^2,$}\\
(\xi\odot\eta)_r=\sin\left(\frac{7}{4}\pi-\frac{1}{2}\sqrt{\frac{9}{4}\pi^2+4x}\right)&\text{ for
$x\in[\pi^2,5\pi^2\slash 2].$}
\end{cases}
\end{equation*}
or in one piece
\begin{equation*}
\xi\odot\eta = \chi_{[0,\frac52\pi^2]}(x)\cdot-\sin(\frac14\pi+\frac14\sqrt{9\pi^2+16x})
\end{equation*}
\begin{figure}[htbp]
\begin{center}
\includegraphics[width=120mm]{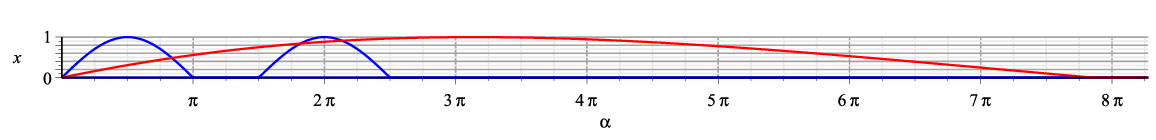}
\caption{$\xi(x)=\chi_{[0,\pi]}(x)\sin x$, $\eta(x)=\chi_{[3\pi\slash 2,5\pi\slash
2]}\cos x,$ and the product $\xi\odot\eta$.}\label{sincos}
\end{center}
\end{figure}

\end{example}

\newpage

\section{$\Ft$}
In section~1 we defined a class of fuzzy numbers $\mathcal{F}_c^0(\R^+)$, where the subscript~``$c$''
stands for ``compact'' to indicate that the characterizing functions defining the fuzzy numbers have compact support,
and the superscript~``$0$'' indicates that these functions are also continuous. \par\noindent
Notably the real numbers $\R$ are beyond the scope of \textit{Definition~1}. In this section a second, wider definition is brought in to include the standard real numbers and also the interval numbers (that is, the characteristic functions of intervals).
Thus, we will allow  a characterizing function $\xi(\cdot)$ to be equal to~$1$ on a closed interval
$[\underline{m},\overline{m}]$ (in the class $\mathcal{F}_c^0(\R^+)$ we had $\underline{m}=\overline{m}$). Our findings and statements will henceforth be formulated for ``fuzzy intervals'' and we will refer to fuzzy numbers as to the subcase when \mbox{$\underline{m}=\overline{m}=m$.}\par
\smallskip
Moreover, we shall no longer require that the fuzzy endpoints $\xi_l$ and $\xi_r$ be strictly monotone and
continuous, and instead require only simple monotonicity and semi-continuity.\par
\smallskip
This new class shall be denoted by $\Ft.$\par

\begin{definition}
A fuzzy interval $\xi$ from the class $\Ft$ is characterized by a function $\xi(\cdot):\R^+\to[0,1]$ with the following
properties:
\begin{itemize}\label{def2}
\item[(1)] compact support: the closure of the points $x\in\R^+$, where $\xi(x)>0$ is a closed interval $[l,r]\subset\R^+,$
\item[(2)] The set of points $x\in\R^+$ where the characterizing function $\xi(\cdot)$ attains value $1$ is a closed interval $[\underline{m},\overline{m}]\subset[l,r]$, in other words\par\noindent  $\{x\in\R:\xi(x)=1\}=[\underline{m},\overline{m}]\subseteq[l,r]$
\item[(3)] the function  $\xi_l(x):= \xi(x)$ for $x\in[l, \underline{m}]$ is non-decreasing and
right-continuous,
\item[(4)] the function  $\xi_r(x):= \xi(x)$ for $x\in[ \overline{m},r]$ is non-increasing and left-continuous.
\end{itemize}
\end{definition}
\begin{remark}
Note that conditions $(1)-(4)$ in Definition~\ref{def2} imply that $\xi(\cdot)$ is upper semi-continuous, that is
$$\limsup_{x\to x_0}\xi(x)\leq\xi(x_0), \quad \forall x_0\in\R^+.$$
\end{remark}\par\noindent
As in section~1 and \textit{Definition 1'} (\ref{lrform1}) the characterizing function of a fuzzy interval $\xi\in\Ft$ can be rewritten as
\par
\smallskip
\noindent{\textbf{Definition 2'}.}
\begin{equation}\label{lrform2}
\xi(x)=\begin{cases}\xi_l(x)&\text{ for $x\in[l,\underline{m}],$}\\
1&\text{ for $x\in[\underline{m},\overline{m}],$}\\
\xi_r(x)&\text{ for $x\in[\overline{m},r].$}
\end{cases}
\end{equation}\par

(Again $\xi_l(\underline{m})=\xi_r(\overline{m})=1$ is an overkill, but we need this for inversion)\par\smallskip

As in the preceding section a fuzzy interval (number) is uniquely determined by its fuzzy components
$\xi_l$ and $\xi_r$. The interval $[\underline{m},\overline{m}]$ need not be stated explicitly, but usually will be for sake of clarity.\par
\smallskip
\textbf{Definition 2''}\label{pair''}
Again a representation by an ordered pair $\left(\xi_l,\xi_r\right)$ where $\xi_l$ and $\xi_r$ meet the conditions of \textit{Definition~2} above seems intuitive and highly desirable, especially when speaking of fuzzy intervals (not fuzzy numbers).\par
\smallskip
As we mentioned before, the class $\mathcal{F}_c(\R^+)$ includes the standard positive numbers $\R^+$ and the
interval numbers. Indeed, every $\lambda\in\R^+$ can be identified with the delta function $\delta_\lambda:\R^+\to\{0,1\}$ such
that $\delta_\lambda(x)=1$ iff $x=\lambda$ and $0$ elsewhere. (\textit{Convention:} We use $\chi_{\{\lambda\}}$ and $\delta_\lambda(x)$ interchangeably). Clearly,
$\delta_\lambda\in\mathcal{F}_c(\R^+)$, which can be seen by
taking $l=\underline{m}=\overline{m}=r=\lambda$ and $\xi_l=\xi_r=\chi_{\{\lambda\}}$ in \textit{Definition~2'}. For an interval number interpreted as its own characteristic function $\chi_{[a,b]}$ of the interval $[a,b]\subset\R^+$ the understanding is the same. We have
$\chi_{[a,b]} \in\mathcal{F}_c(\R^+)$, which can be seen by taking $l=\underline{m}=a$ and $b=\overline{m}=r$ and $\xi_l=\chi_{\{a\}},\xi_r=\chi_{\{b\}}$
in \textit{Definition~2'}.

\smallskip
We shall see later that the arithmetic which will be defined below for $\Ft$ is consistent with the standard arithmetic of
$\R^+$, as well as the set arithmetic for intervals.\par

\begin{subsection}{Generalized Inverse Function}\label{GIF}
The (not strictly) monotone and/or possibly discontinuous functions $\xi_l$, $\xi_r$ of \textit{Definition~2'} may not be invertible in the classical sense. In order to be able to employ, as before in section~1 formulae analogous to \eqref{dodawanie} and \eqref{mnozenie} in defining the operations ``$\oplus$'',``$\odot$'' we shall presently introduce and then apply the notion (appearing in statistics in the form of generalized inverses of CDFs, or quantile functions) of a \emph{generalized inverse function}:\par
\smallskip\noindent
Take a function $h:I\mapsto[0,1]$ defined on a closed interval $I$, that is one-sided continuous and non-strictly
monotone.par\noindent
The price to pay for relaxing the conditions of continuity and strict monotonicity and instead assuming only one-sided continuity and non-strict
monotonicity is that we encounter two types of problems in the process of finding a useful inverse function in the classical sense:\par
\noindent
\begin{enumerate}
\item
\vspace{1mm} There is a point $x_0$, where the function $h$ has a jump discontinuity.\vspace{1mm}

\noindent Such a function function is invertible in the classical sense, but the domain of its inverse function will not be the full interval $[0,1]$.\vspace{0.5mm}

\noindent Now two functions ${h_1}, {h_2}$ with jump discontinuities at different points will produce inverse functions ${h_1}^{-1}, {h_2}^{-1}$ of different domains~$h_1(I)$ and $h_2(I).$ -
Therefore, \eqref{dodawanie}, and \eqref{mnozenie} cannot be applied because the addition or multiplication of functions with different domains is  simply not well defined.
\item
The function $h$ is constant on some interval $[\underline{x}_0,\overline{x}_0]$.\vspace{1mm}

\noindent In this case the function is not injective, and therefore not invertible at all in the classical sense.
\end{enumerate}\smallskip

In order to overcome these problems and to be able to extend formulae \eqref{dodawanie} and \eqref{mnozenie} we shall, for our present purposes, single out two types of generalized inverse functions. For a non-decreasing right-continuous function $f:I\to J$ defined on
on a closed interval $I$ and taking values in a closed interval $J$ (resp. non-increasing left-continuous function
$g:I\to J$) we define two types of \emph{generalized inverse functions} by

\begin{equation}\label{geninvinf}
f_{\inf}^{-1}(\alpha) =
\begin{cases}
\vspace{1mm}
\inf\{x: \xi_{l}(x) \geq \alpha\} \quad\text{for}\, \alpha \in (0,1],\\
\underline{\supp(\xi_{l})}\quad\text{for}\, \alpha = 0,
\end{cases}
\end{equation}
and
\begin{equation}\label{geninvsup}
f_{\sup}^{-1}(\alpha) =
\begin{cases}
\vspace{1mm}
\sup\{x: \xi_{r}(x) \geq \alpha\} \quad\text{for}\, \alpha \in (0,1],\\
\overline{\supp(\xi_r)}\quad\text{for}\, \alpha = 0.
\end{cases}
\end{equation}
\smallskip

To see how concretely the generalized inverse procedure works in practice let us consider the case of a non-decreasing
right-continuous function $f$ (which may serve as the left fuzzy endpoint $\xi_l$ of a fuzzy interval $\xi$ as in
\textit{Definition~\ref{def2}}).

\begin{enumerate}
\item
Let $f$ be a non-decreasing right-continuous function with a jump discontinuity at the point $x_0$. Let
$f(x_0)=\alpha_0$. Let us set $\overline{\alpha}_{0}=\alpha_0$ and $\underline{\alpha}_0=\sup\{f(x):x<x_0\}$.
Then for all $\alpha\in(\underline{\alpha}_0,\overline{\alpha}_0]$ we have that
$$f_{\inf}^{-1}(\alpha)=\inf\{x:f(x)\geq\alpha\}=x_0,$$ that is, the generalized inverse function is constant on
the half-closed interval  $(\underline{\alpha}_0,\overline{\alpha}_0]$.\vspace{1.2mm}

\item
Let $f$ be constant on the interval $[\underline{x}_0,\overline{x}_0)$ with $f(x)=\alpha_0$ for $x$ in the
interval. In this case clearly,
 $$f_{\inf}^{-1}(\alpha_0) =\inf\{x:f(x)\geq\alpha_0\}= \underline{x}_0,$$
 while
 $$
 \inf\{f_{\inf}^{-1}(\alpha):\alpha>\alpha_0\}=\overline{x}_0.
 $$
 Thus, the generalized inverse has a jump discontinuity at $\alpha_0$.
\end{enumerate}

\noindent For a (not strictly) decreasing function $g$, which in what is coming may serve as the right fuzzy component of a fuzzy interval $\xi$ the reasoning is analogous, except that ``$\inf$" on each occurrence must be replaced by ``$\sup$", and right-closed intervals be substituted by left-closed ones and vice-versa.\par\smallskip

In the light of the above we see that the procedure of finding the generalized inverse works may be unceremoniously summarized: Jump discontinuities of the input convert to intervals where the output inverted function is constant. And the intervals where
the input is constant convert to jump discontinuities of the output. These two observations are illustrated below for fragments of two possible left fuzzy endpoints (green) and their inverses (blue) :
\begin{figure}[H]
\includegraphics[width=60mm]{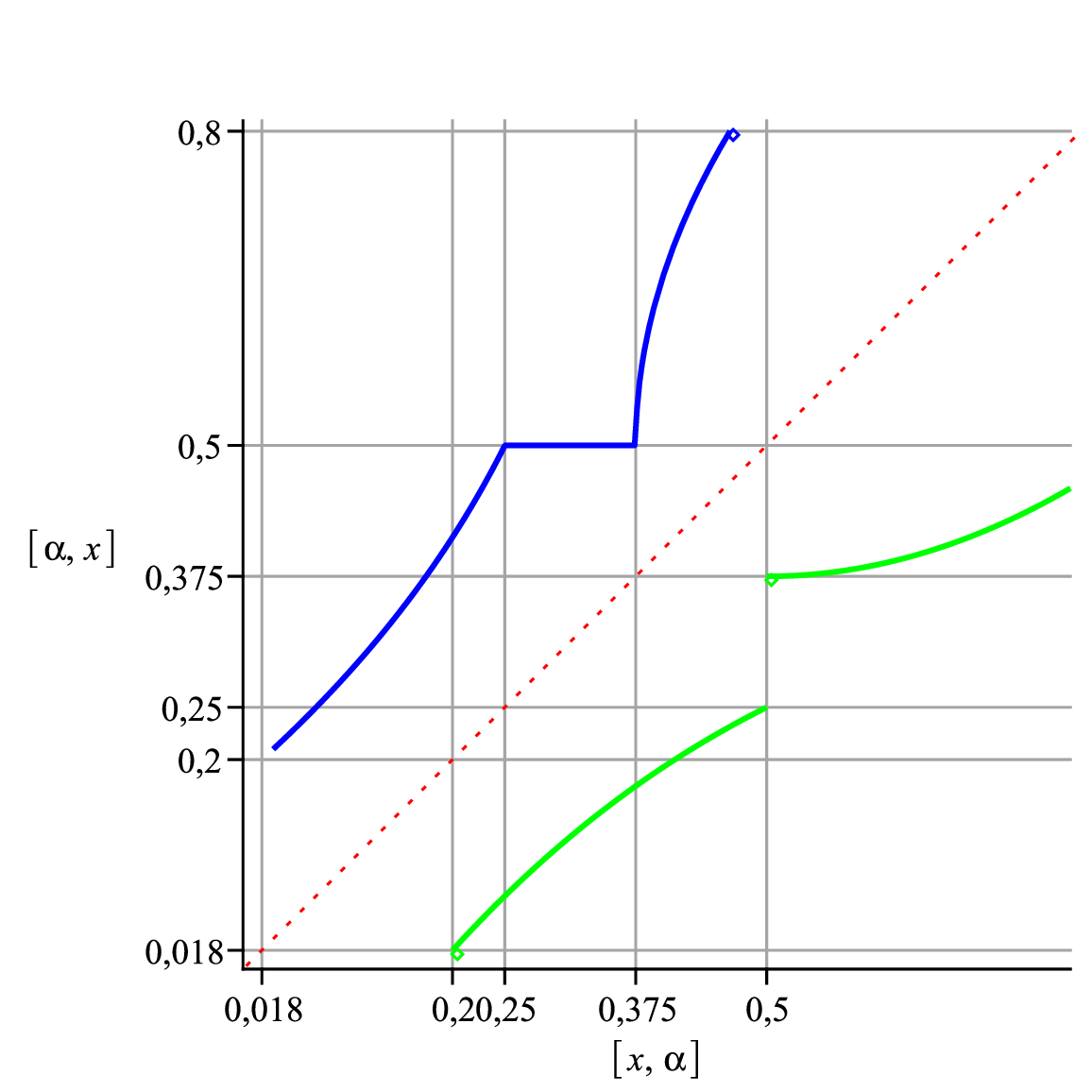}
\includegraphics[width=60mm]{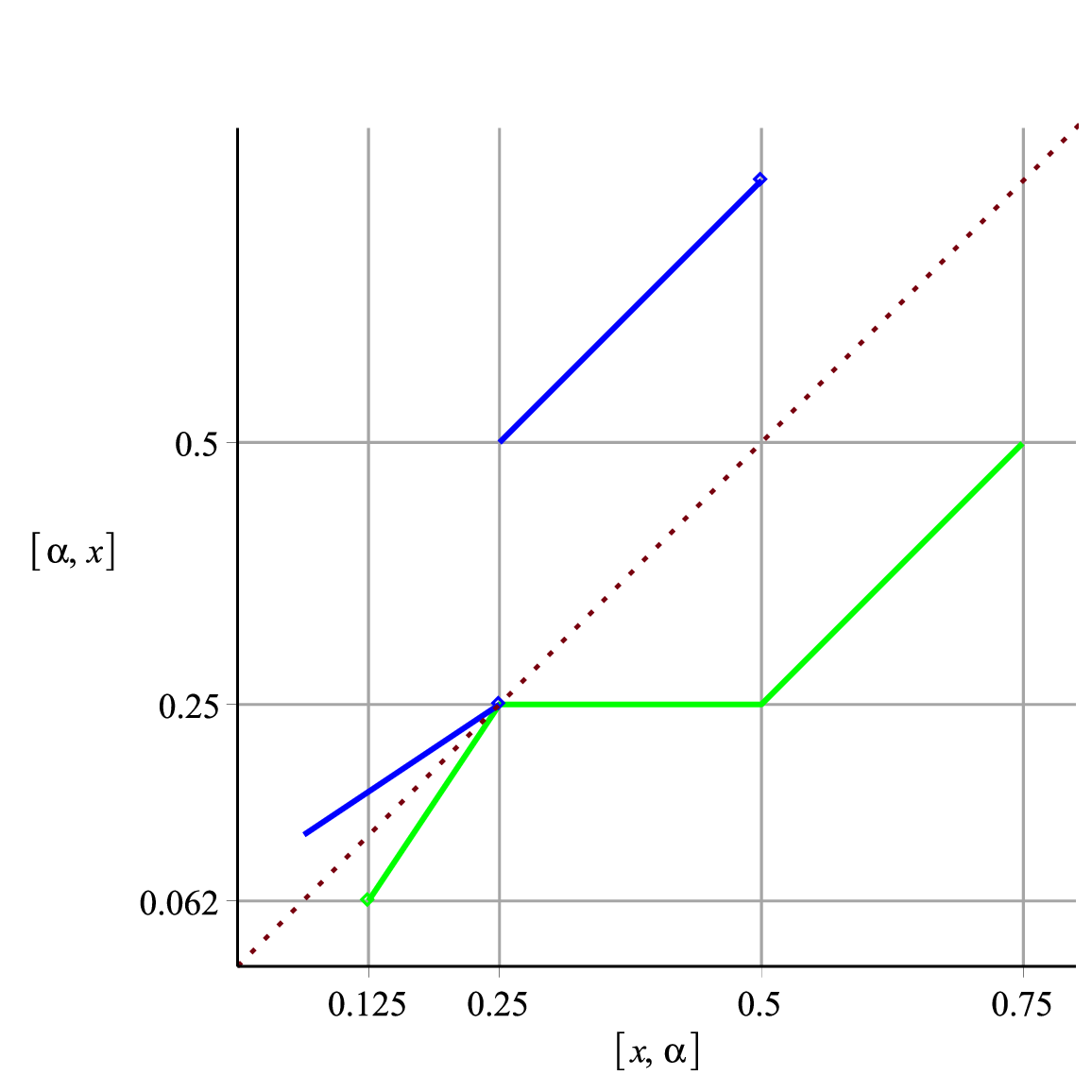}
\caption{\tiny jump-discontinuity converts to constant, constant to jump-discontinuity}
\end{figure}

Again for a non-increasing left-continuous function $g$ serving as the right fuzzy part $\xi_r$ of a fuzzy interval
$\xi$ given as in \textit{Definition~2} the procedure of finding the general inverse works analogously substituting ``$\sup$" for ``$\inf$".\par

\smallskip
\begin{remark}\label{invsemicont}
We also explicitly note, that for a left fuzzy component $\xi_l(\cdot)$ which is by definition increasing and right-continuous its generalized
inverse\par\noindent ${\xi_l}^{-1}_{\inf}:[0,1]\to\R^+$ is still increasing (only fuzzy convexity i.e.: $\xi(\lambda x_1+(1-\lambda)x_2)\geq\min\bigl(\xi(x_1),\xi(x_2)\bigr)$ transforms into concavity), but left-continuous.\par
\noindent Observe also that ${\xi_l}^{-1}_{\inf}(0)=l$, ${\xi_l}^{-1}_{\inf}(1)=\underline{m}$.\par\noindent
\vspace{1mm}

\noindent Likewise ${\xi_r}^{-1}_{\sup}:[0,1]\to\R^+$ is a decreasing, right-continuous function and\par
\noindent ${\xi_r}^{-1}_{\sup}(1)=\overline{m}$ and ${\xi_r}^{-1}_{\sup}(r)=0.$
\end{remark}

\bigskip

\begin{subsubsection}{Graphs illustrating generalized inversion}$ $\par

\smallskip
In what follows some example graphs are displayed:
\newpage

\noindent The infimum and supremum inverse of a real number $\lambda$ understood as the fuzzy number $\chi_{\{\lambda\}}$ or equivalently the delta function $\delta_\lambda$ is the constant function $\alpha=\lambda$  (with domain $[0,1]$):
\begin{figure}[H]
        \centering

                \centering
                \includegraphics[width=0.6\textwidth]{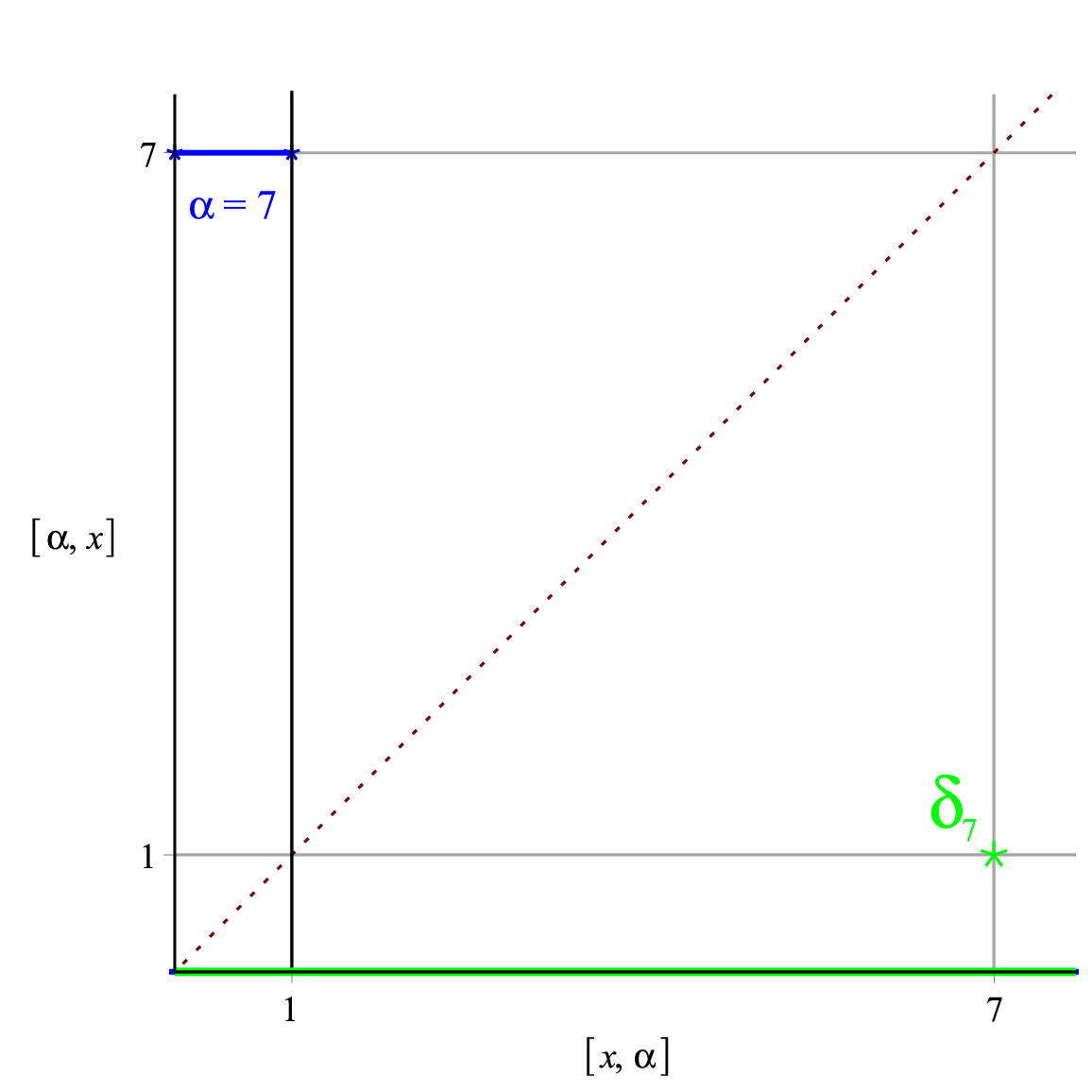}
                \caption{inverse (infimum and supremum) of a real number}
                \label{fig:delta}

\end{figure}
\bigskip
\noindent Below the infimum and supremum generalized inverses (blue) of a stepfunction (green) are confronted for comparison:
\begin{figure}[H]
\begin{subfigure}[b]{0.5\textwidth}
                \centering
                \includegraphics[width=\textwidth]{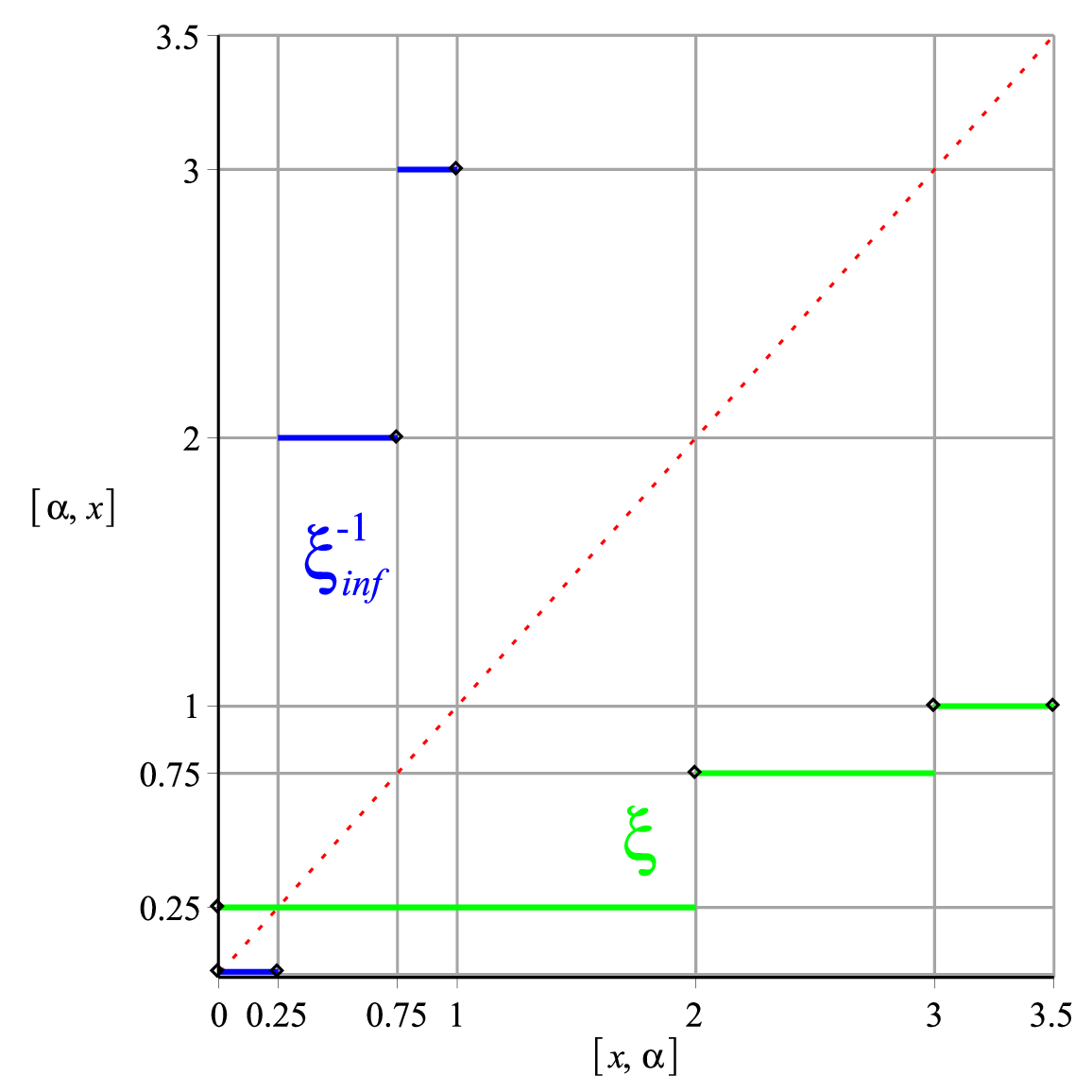}
                \caption{\tiny infimum inverse of a step function (green)}
                \label{fig:stepfunction1}
        \end{subfigure}%
\begin{subfigure}[b]{0.5\textwidth}
                \centering
                \includegraphics[width=
                \textwidth]{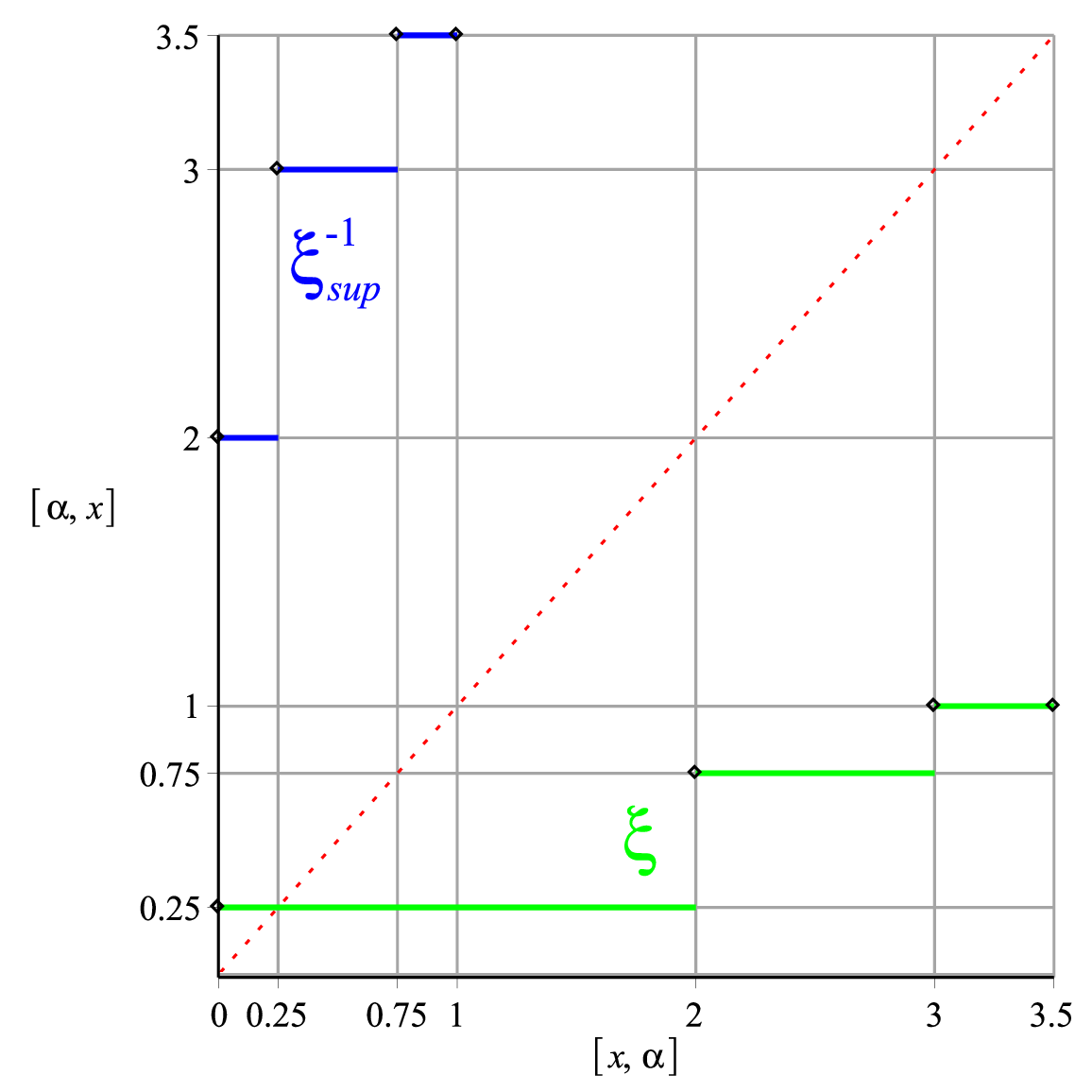}
                \caption{\tiny supremum inverse of the same step function}
                \label{fig:stepfunction2}
\end{subfigure}
\end{figure}

\newpage

As a final example illustrating the procedure of generalized inversion here is now a freely patched together ``fantasy" function which
includes both jumps and constant parts and its generalized inverse:\par

\bigskip
\begin{equation*}
{fanta}_{l}(x)=\begin{cases}
\frac14 x\quad &x\in[\frac{1}{4},\frac23)\\ \sin(3x-\frac12) \quad &x\in[\frac23,\frac67) \\ \frac58 \quad &x\in[\frac67,\frac98) \\
\frac56x \quad &x\in[\frac98,\frac65]
\end{cases}\end{equation*}\par\noindent

\begin{figure}[H]
 \centering
\includegraphics[width=0.6\textwidth]{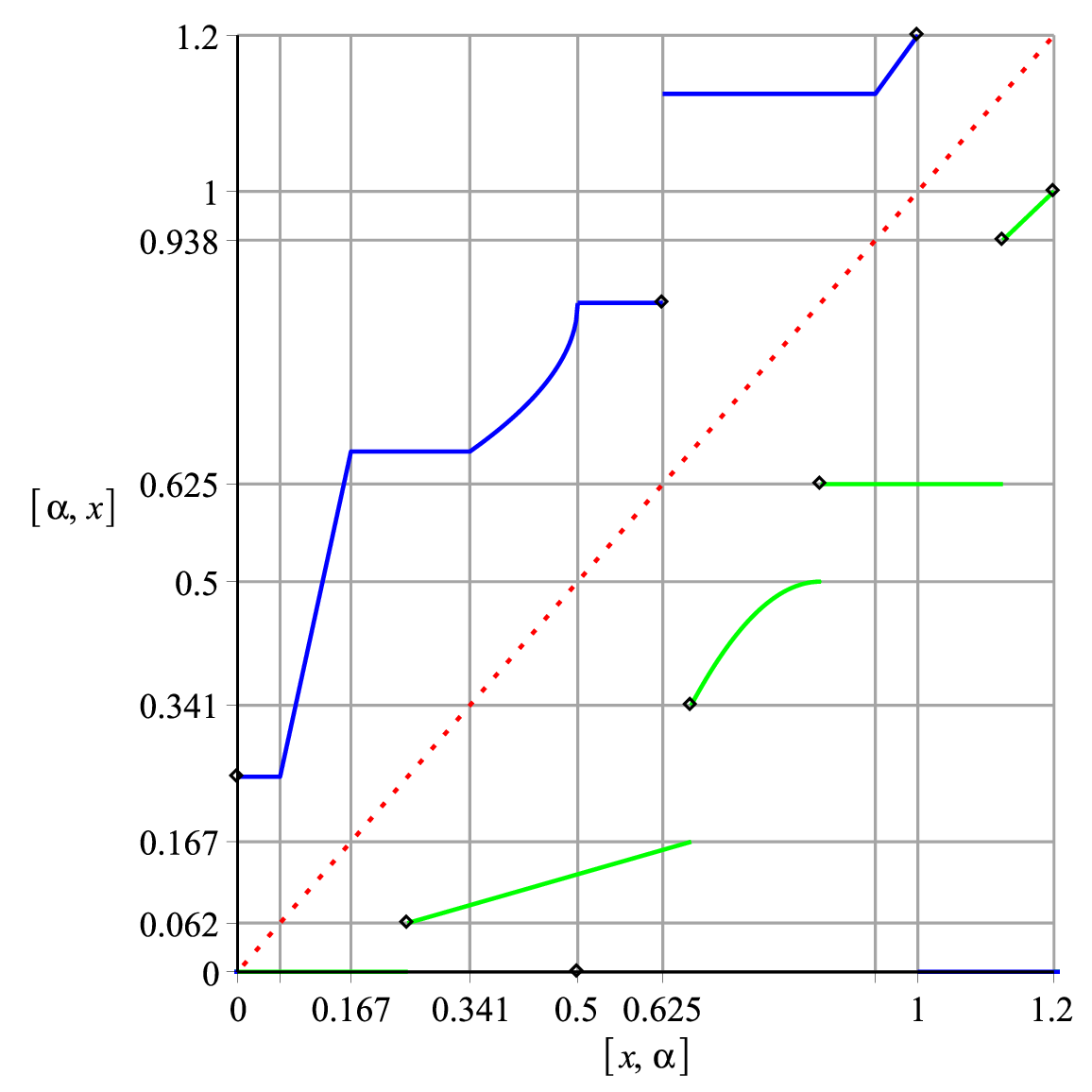}
\caption{the infimum generalized inverse (blue) of a ``fantasy" function (green)}
\end{figure}

\begin{equation*}
{fanta}_{d}(\alpha)=\begin{cases}
\frac14 \quad &\alpha\in [0,\frac{1}{16}]\\
4\alpha \quad &\alpha\in (\frac1{16},\frac16]\\ \frac23 \quad &\alpha\in(\frac16,\sin(1)-\frac12] \\ \begin{tiny}\frac13 \arcsin(\alpha+\frac12)+\frac13 \end{tiny}
&\alpha\in(\sin(1)-\frac12,\sin(\frac{11}7)-\frac12] \\ \frac16 \quad &\alpha\in(\sin(\frac{11}{7})-\frac12,\frac58] \\
\frac98 \quad &\alpha\in(\frac58,\frac{15}{16}] \\ \frac65\alpha \quad &\alpha\in(\frac{15}{16},1]
\end{cases}\end{equation*}

\end{subsubsection}
\end{subsection}
\newpage

\subsection{Re-inversion of fuzzy endpoints}\label{opardef2}$ $\par\vspace{1mm}
\noindent To reduce the multitude of sub- and superscripts and indices let us first introduce some shorthand notation:

\begin{notation}$ $\par

\begin{enumerate}
\item Set ${\xi_{l}}^{-1}_{\inf} =:\xi_d$ as in \textit{downwards} fuzzy, and
\item set ${\xi_{l}}^{-1}_{\inf} =:\xi_d$ as in \textit{upwards} fuzzy
\end{enumerate}

\end{notation}

To be able to adapt formulae~\ref{dodawanie} and \ref{mnozenie} to the present case we need to define an suitable re-inversion-operator for the generalized inversion function. One basic prerequisite condition which must be expected from this operator, which we will denote by and overhead left arrow $``\overleftarrow{(\cdot)}"$ is that $$\overleftarrow{(\xi_d)}(x) = \xi_{l}(x) \text{\quad and \quad} \overleftarrow{(\xi_u)}(x) = \xi_{r}(x).$$

\begin{subsubsection}{Inversion and re-inversion of a fuzzy left endpoint $\xi_l(\cdot)$}$ $\par\noindent

\smallskip\noindent
For easy reference recall that for $\xi_l: \R\mapsto[0,1]$ onto, non-decreasing and left-continuous
its generalized inverse function is defined by:\vspace{2mm}
\begin{equation}\label{geninvinf'}
\xi_d(\alpha) =
\begin{cases}
\vspace{1mm}
\inf\{x: \xi_{l}(x) \geq \alpha\} \quad\text{for}\: \alpha \in (0,1],\\
\underline{\supp(\xi_{l})}\quad\text{for}\: \alpha = 0.
\end{cases}
\end{equation}
$ $\par

\smallskip
\noindent
\begin{remark}
Note, that for $\xi_l$ onto, strictly monotonic (increasing) and continuous the definitions of generalized and classical inversion coincide:
\begin{equation}
\xi_d = \xi_l^{-1}.
\end{equation}
\end{remark}

\begin{remark}
Remember, that the generalized inverse $\xi_d$ of the non-decreasing, left-continuous function $\xi_l$ is again a non-decreasing function, but right-continuous.\vspace{-3.mm}
\end{remark}

\bigskip\noindent
Now define a re-inversion operator for the generalized inverse of a left fuzzy endpoint by: \vspace{2mm}
\begin{equation}
\overleftarrow{\xi_d}(x) =
\begin{cases}
\vspace{1mm}
\inf\{\alpha: \xi_{d}(\alpha)>x\} \quad\text{for}\, x\in\left[\underline{\supp(\xi_l)}, \overline{\supp(\xi_l)}\right),\\
\vspace{1mm}
1  \quad\text{for}\, x = \overline{\supp(\xi_l)},\\
0  \quad\text{for}\, x \notin \supp(\xi_l).
\end{cases}
\end{equation}

\noindent It is then easy to see that
\begin{lemma}
\begin{equation}
\overleftarrow{(\xi_d)}(x) = \xi_{l}(x).\label{maintheorem}
\end{equation}
\end{lemma}

\begin{proof}
$ $  The proof follows the reasoning applied in the preceding section~\ref{GIF}
\end{proof}

\begin{remark}
Note that for the strictly monotonic and continuous case \eqref{maintheorem} reads just
\begin{equation}
(\xi_l^{-1})^{-1} = \xi_l.
\end{equation}
\end{remark}

\end{subsubsection}\vspace{2mm}

\noindent We now proceed analogously with the right fuzzy endpoint:\par\vspace{1mm}\noindent

\begin{subsubsection}{Inversion and re-inversion of a right fuzzy endpoint $\xi_r(\cdot)$}$ $\par\noindent

\noindent We recall that the generalized inverse to the non-increasing, left-continuous right fuzzy endpoint $\xi_r: \R\mapsto[0,1]$ is given by:

\begin{equation}\label{geninvsup'}
\xi_u(\alpha) =
\begin{cases}
\vspace{1mm}
\sup\{x: \xi_{r}(x) \geq \alpha\} \quad\text{for}\, \alpha \in (0,1],\\
\overline{\supp(\xi_r)}\quad\text{for}\, \alpha = 0.
\end{cases}
\end{equation}
$ $
\begin{remark}
Obviously for $\xi_r$ strictly decreasing and continuous the equality
\begin{equation}
\xi_u = \xi_r^{-1}
\end{equation}
will hold.
\end{remark}

\begin{remark}
Note, that the generalized inverse $\xi_u$ of the non-increasing, right-continuous function $\xi_l$ is again a non-increasing function, but left-continuous.\vspace{-3.mm}
\end{remark}

\smallskip
\noindent And we define re-inversion by setting:

\begin{equation}
\overleftarrow{(\xi_u)}(x) =
\begin{cases}
\vspace{1mm}
\sup\{\alpha: \xi_{u}(\alpha)>x\} \quad\text{for}\, x\in\supp(\xi_r),\\
0  \quad\text{for}\, x \notin \supp(\xi_r).
\end{cases}
\end{equation}

\smallskip
\begin{lemma}
Then: $\overleftarrow{(\xi_u)}(x) = \xi_r(x).$
\end{lemma}
\begin{proof}
As above.
\end{proof}

\end{subsubsection}

With this the requisite apparatus for extending the definitions of section~1 has been established and we may proceed:

\begin{subsection}{Addition and multiplication in $\Ft$}$ $

\smallskip
Addition and multiplication in $\Ft$ are defined in the same way as in the class $\Fo$, namely we extend formulae
\eqref{dodawanie} and \eqref{mnozenie} by application of the generalized inverse functions as defined in \eqref{geninvinf'} and \eqref{geninvsup'}:\par
\smallskip We set
\begin{equation}\label{uododawanie}
\begin{split}
(\xi\oplus\eta)_{l}=\overleftarrow{(\xi_d+\eta_d)}\\
(\xi\oplus\eta)_{r}=\overleftarrow{(\xi_u+\eta_u)}
\end{split}
\end{equation}

and

\begin{equation}\label{uomnozenie}
\begin{split}
(\xi\odot\eta)_{l}=\overleftarrow{(\xi_d\cdot\eta_d)}\\
(\xi\odot\eta)_{r}=\overleftarrow{(\xi_u\cdot\eta_u)}
\end{split}
\end{equation}

\end{subsection}

\smallskip
\begin{notation}
To keep the text homogeneous in the further the above notation, i.e. subscripts ``d,u'', for generalized inversion of discontinuous and/or not strictly monotonous fuzzy endpoints, and the overhead arrow symbol ``$\overleftarrow{ }$'' for generalized re-inversion will be applied also to endpoints which are actually invertible in the classical sense.
\end{notation}

\begin{remark} Clearly, the sum and the product of left-continuous (resp. right-continuous) functions are
left-continuous (resp. right-continuous), thus the above operations are well defined.
\end{remark}

\begin{example}[Multiplication of interval numbers]
The basic example of a fuzzy interval is the abstraction of a real interval $A=[\underline{a},\overline{a}]$ by its characterising function $\xi_{A}=\chi_{[\underline{a},\overline{a}]}$.\par
Let us fuzzy multiply two interval numbers $\xi_{A}=\chi_{[2,3]}$ and $\xi_{B}=\chi_{[5,6]}.$ The result will be, as expected,
$\xi_{A}\odot\xi_{B}=\chi_{[10,18]}$, as shown
in Fig.~\ref{odc}.

\noindent Going by the definiton, one piece at a time, step by step:
\begin{equation}
\begin{split}
{\xi_A}_{\,l}=\chi_{\{2\}} \text{ for } x\in\R\\
{\xi_B}_{\,l}=\chi_{\{5\}} \text{ for } x\in\R\\
{\xi_A}_{\,d}(\alpha)=2 \text{ for }\alpha \in [0,1]\\
{\xi_B}_{\,d}(\alpha)=5 \text{ for }\alpha \in [0,1]\\
{\xi_A}_{\,d}(\alpha)\cdot{\xi_B}_{\,d}(\alpha)=2\cdot5=10 \text{ for }\alpha \in [0,1]\\
(\xi_A\odot\xi_B)_{l}=\overleftarrow{{\xi_A}_{\,d}(\alpha)\cdot{\xi_B}_{\,d}(\alpha)}=\overleftarrow{(\alpha=10)}=
\chi_{\{10\}}\text{ for }x \in \mathbb{R}
\end{split}
\end{equation}

The computation of the righthand side $(\xi_{A}\odot\xi_{B})_r$ is analogous.

\begin{figure}[H]
\begin{center}
\includegraphics[height=20mm]{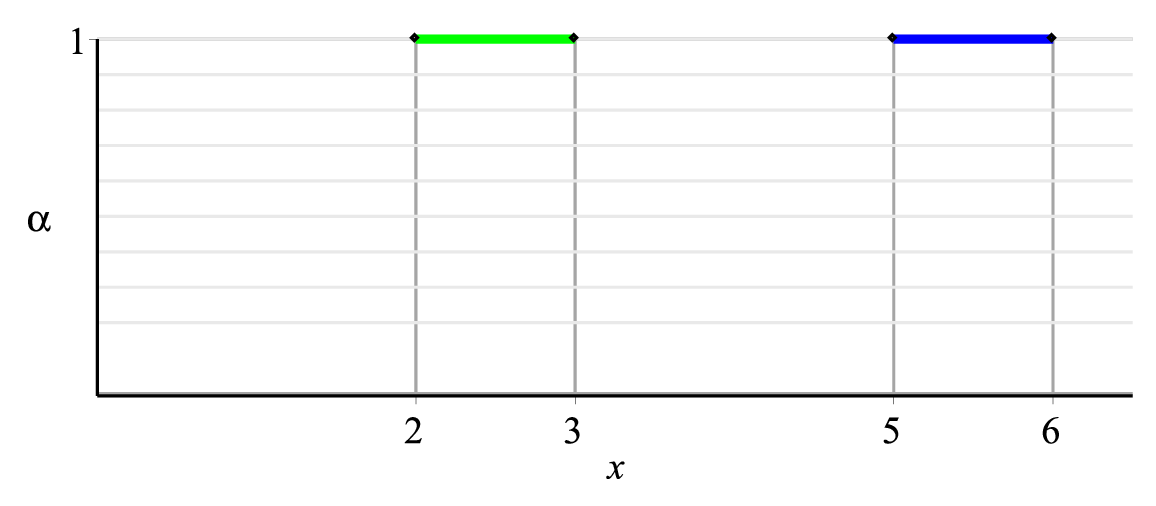}
\includegraphics[width=20mm]{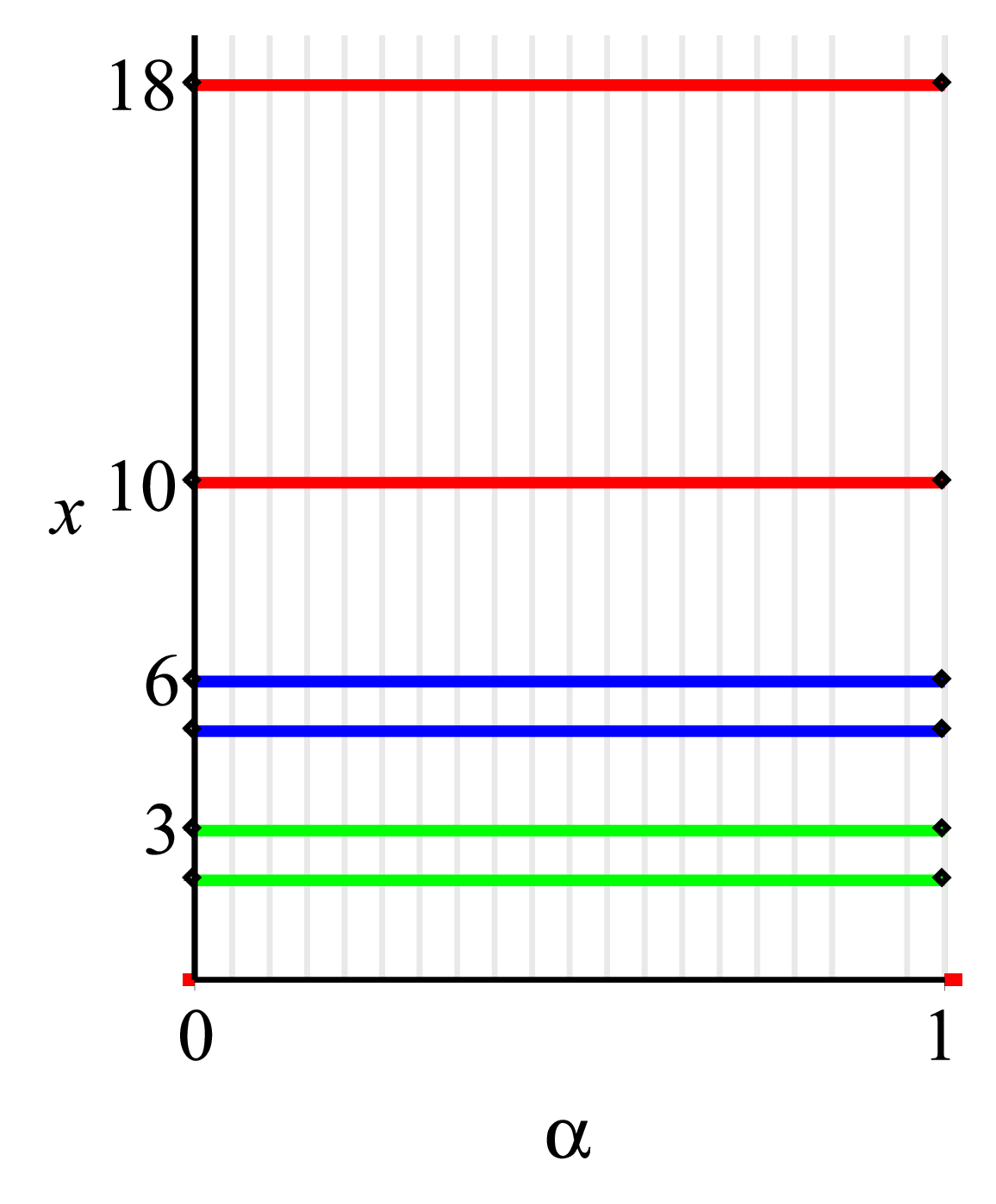}
\includegraphics[height=20mm]{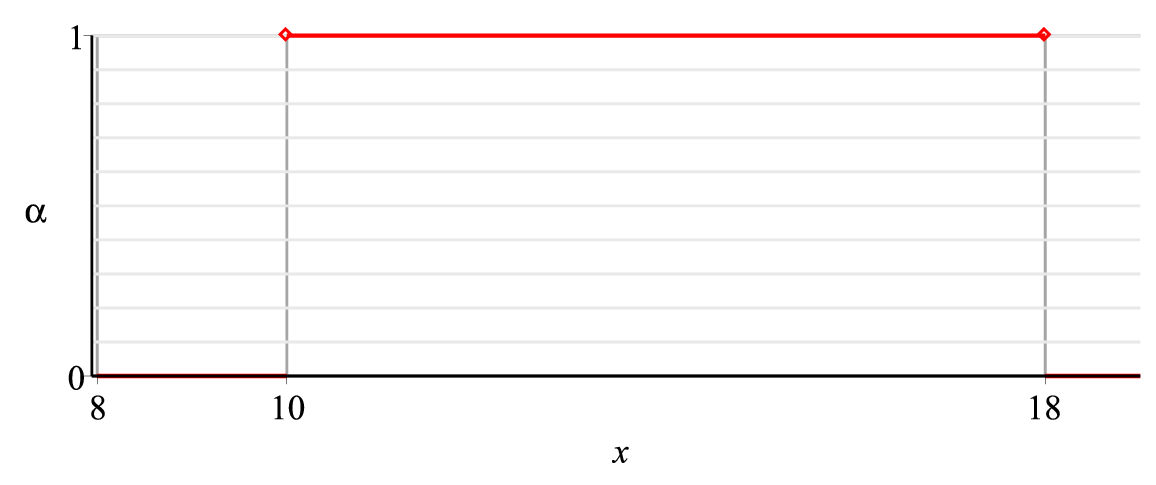}
\caption{\tiny Multiplication of interval numbers: (a) $\xi_{A}=\chi_{[2,3]}$ and $\xi_{B}=\chi_{[5,6]},$ (b) the inverted intervals and their product, (c) the result $\xi_{A}\odot\xi_{B}=\chi_{A\cdot B}=\chi_{[10,18]}$.}\label{odc}
\end{center}
\end{figure}

\smallskip\noindent
In general, for $0\leq\underline{a}\leq\overline{a}$ and  $0\leq\underline{b}\leq\overline{b}$ we may always use
\begin{equation}\label{przedzialowe}
\begin{split}
\chi_{[\underline{a},\overline{a}]}\odot\chi_{[\underline{b},\overline{b}]}=&\chi_{[\underline{a}\underline{b},\overline{a}\overline{b}]},\\
\chi_{[\underline{a},\overline{a}]}\oplus\chi_{[\underline{b},\overline{b}]}=&\chi_{[\underline{a}+\underline{b},\overline{a}+\overline{b}]}.
\end{split}
\end{equation}
which is of course perfectly consistent with classical interval arithmetic on the nonnegative reals.
\end{example}\par

Before proceeding to multiplication by a real number let us introduce some simplifying shorthand:$ $\vspace{1mm}

\noindent{\textit{Notation:}} \normalsize We write ${\delta_{\lambda}}^{-1}_{\inf} = {\delta_{\lambda}}_d = {\delta_{\lambda}}^{-1}_{\sup} = {\delta_{\lambda}}_u = \lambda \text{ for } \alpha\in[0,1]$  or simply $\alpha=\lambda$\par
\smallskip

\begin{example}[Multiplication of a fuzzy interval by a real number]
We calculate for $x\in\bigl[\underline{\supp{(\lambda\odot\xi)_l}}, \overline{\supp{(\lambda\odot\xi)_l}}\bigr)$
\begin{equation*}
\begin{split}
(\lambda\odot\xi)_l(x)=& \; \overleftarrow{({\delta_{\lambda}}_{\,d}\cdot{\xi_d})}(x)=\overleftarrow{(\lambda\cdot{\xi_d}})(x)\\
= & \quad\inf\{\alpha:(\lambda\cdot\xi_d)(\alpha)> x\}\\
= & \quad\inf\{\alpha:\xi_d(\alpha) > x\slash\lambda\}\\
= & \quad\overleftarrow{\xi_d}(x\slash\lambda)\\
= & \quad\xi_l(x\slash\lambda).
\end{split}
\end{equation*}

The case $x = \overline{\supp{(\lambda\odot\xi)_l}}$ is trivial.
Repeating the calculation for $\xi_r$, we have established the following property:
\begin{equation}\label{mnozenie przez rzeczywista}
(\lambda\odot\xi)(x)=\xi(x\slash\lambda).
\end{equation}\par

\smallskip\noindent
Below $$2\odot\chi_{[3,4]}(x)\sqrt{4-x} =\chi_{[6,8]}(x)\sqrt{4-x/2}$$ is shown as an example:
\begin{figure}[H]
\begin{center}
\includegraphics[width=70mm]{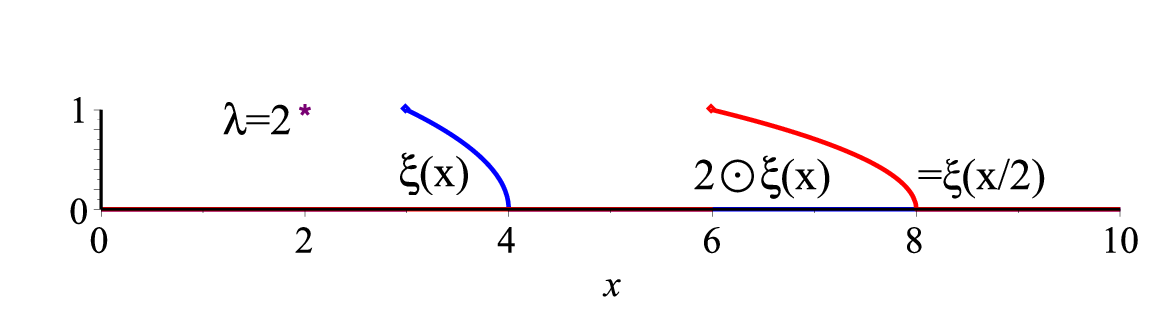}
\caption{\tiny Multiplication of a fuzzy interval by a real number}\label{po}
\end{center}
\end{figure}

\end{example}

\noindent Property \eqref{mnozenie przez rzeczywista} provides for the following formula which is abecedarian for statistical
analysis of fuzzy data:

\begin{equation}
\frac{1}{n}\odot(\xi\oplus\xi\oplus\dots\oplus \xi)=\frac{1}{n}\odot\bigoplus_{i=1}^n\xi=\xi,
\end{equation}

Indeed,
\begin{equation*}
\frac{1}{n}\odot\bigoplus_{i=1}^n\xi_r(x)=\frac{1}{n}\odot\overleftarrow{\left(n\odot\xi_u(x)\right)}=\frac{1}{n}\odot{\xi_r}(x\slash
n)=\xi_r(x),
\end{equation*}
and a similar check can be done for $\xi_l(x)\,.$

\begin{example}[Adding a real number to a fuzzy interval]
We have seen above that multiplying a fuzzy interval by a real number yields an inverse dilation of the fuzzy interval. Similarly, it can be easily
shown that adding a real number to a fuzzy interval corresponds to a negative translation:
\begin{equation}\label{dodawanie rzeczywistej}
(\lambda\oplus\xi)(x)=\xi(x-\lambda)\:.
\end{equation}
\end{example}

\begin{example}
In the last example of this subsection we shall provide calculations and graphs of addition and multiplication of
a fuzzy interval $\xi$ and a fuzzy number $\eta$.\vspace{1mm}

For $\eta$ we take the already considered function $\chi_{[0,\pi]}\cdot\sin$:

\begin{equation}
 \eta = \begin{cases}
         \eta_l=\sin(x)\quad &x\in[0,\frac{\pi}{2}]\\
         1&x=\frac{\pi}{2}\\
         \eta_r=\sin(x) \quad &x\in[\frac{\pi}{2},\pi]\:.
         \end{cases}
 \end{equation}\vspace{1mm}

 This example's $\xi$ has a linear part $\xi_l$ and a step function as $\xi_r$\,:

 \begin{equation}
 \xi = \begin{cases}
         \xi_l=8x-1 \quad& x\in[\frac{1}{8},\frac{1}{4}]\\\vspace{1mm}
        1\quad &x\in[\frac{1}{4},1]\\
         \xi_r=\tiny
           \begin{cases}
           1&x=1\\
           0.3&x\in(1,2]\\
           0.1&x\in(2,4] \,.
           \end{cases}
        \end{cases}
 \end{equation}

\noindent Here the function $\xi_r$ is not invertible in the classical sense.\smallskip

Both characterizing functions are pictured below in one graph. The middle parts are marked in dotted grey and the left and right fuzzy components are respectively green and blue.

\begin{figure}[H]
\begin{center}
\includegraphics[width=100mm]{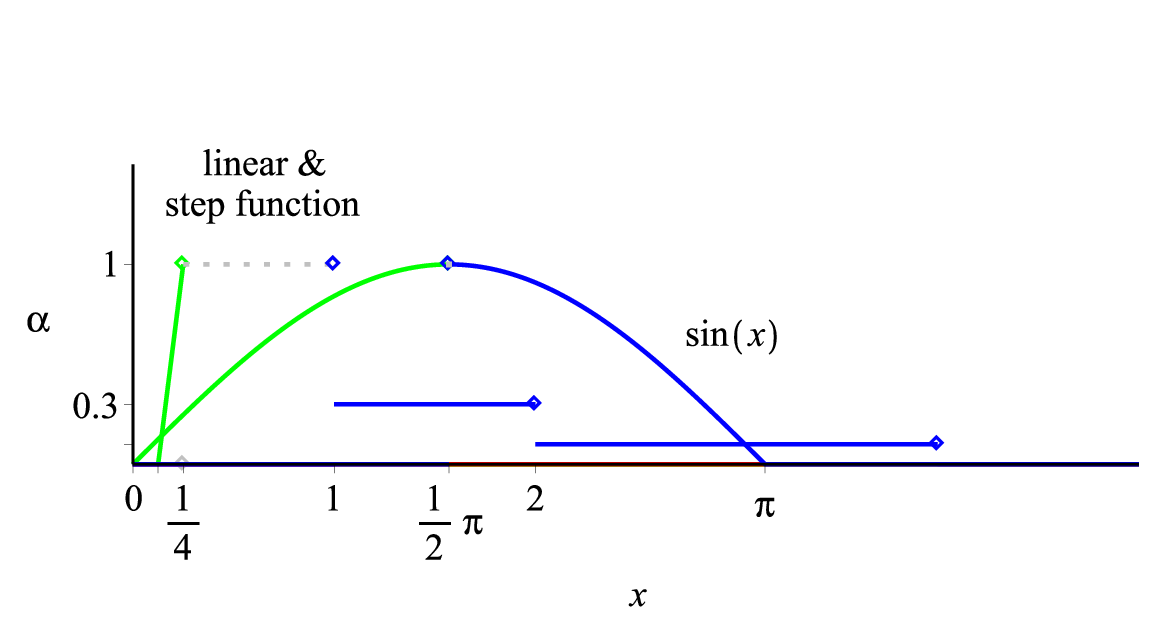}
\end{center}
\caption{sine, linear and step function}
\end{figure}

We have
\begin{align*}
{\xi_d}(\alpha)= \frac{1-\alpha}{8} \qquad \text{ \text{for} $\alpha\in[0,1]$ }
\end{align*}
and
\begin{equation*}
{\xi_u}(\alpha)=
\begin{cases}
4 \quad\text{  \text{for} } \alpha\in[0,0.1)\\
2 \quad\text{  \text{for} } \alpha\in[0.1,0.3)\\
1 \quad\text{  \text{for} } \alpha\in[0.3,1]\,.
\end{cases}
\end{equation*}
\par

\smallskip
\noindent Let us write out  $\eta^{-1}(\alpha)$ as in Example~\ref{sinus}:
\begin{equation*}
\begin{split}
{\eta_d}(\alpha)=& \arcsin(\alpha) \quad\text{ for $\alpha\in[0,1],$}\\
{\eta_u}(\alpha)=& \pi-\arcsin(\alpha) \quad\text{ for $\alpha\in[0,1]$.}
\end{split}
\end{equation*}

We receive
\begin{equation*}
\xi_d+\eta_d=\frac{1-\alpha}{8}+\arcsin(\alpha) \qquad \text{ for $\alpha\in[0,1]$ }
\end{equation*}
and
\begin{equation*}
\xi_u+\eta_u=
\begin{cases}
4 + \pi-\arcsin(\alpha)\quad\text{  \text{for} } x\in[0,0.1)\\
2 + \pi-\arcsin(\alpha)\quad\text{  \text{for} } x\in[0.1,0.3)\\
1 + \pi-\arcsin(\alpha)\quad\text{  \text{for} } x\in[0.3,1].
\end{cases}
\end{equation*}
\bigskip

Multiplication:
\begin{equation*}
\xi_d\cdot\eta_d=\frac{1-x}{8}\cdot\arcsin(\alpha) \qquad \text{ for $\alpha\in[0,1]$}
\end{equation*}
and
\begin{equation*}
\xi_u\cdot\eta_u=\begin{cases}
4\cdot(\pi-\arcsin(\alpha))\quad\text{ for } \alpha\in[0,0.1)\\
2\cdot(\pi-\arcsin(\alpha))\quad\text{ for } \alpha\in[0.1,0.3)\\
1\cdot(\pi-\arcsin(\alpha))\quad\text{ for } \alpha\in[0.3,1]\,.
\end{cases}
\end{equation*}
\bigskip

\begin{figure}[H]
\begin{center}
\includegraphics[width=23mm]{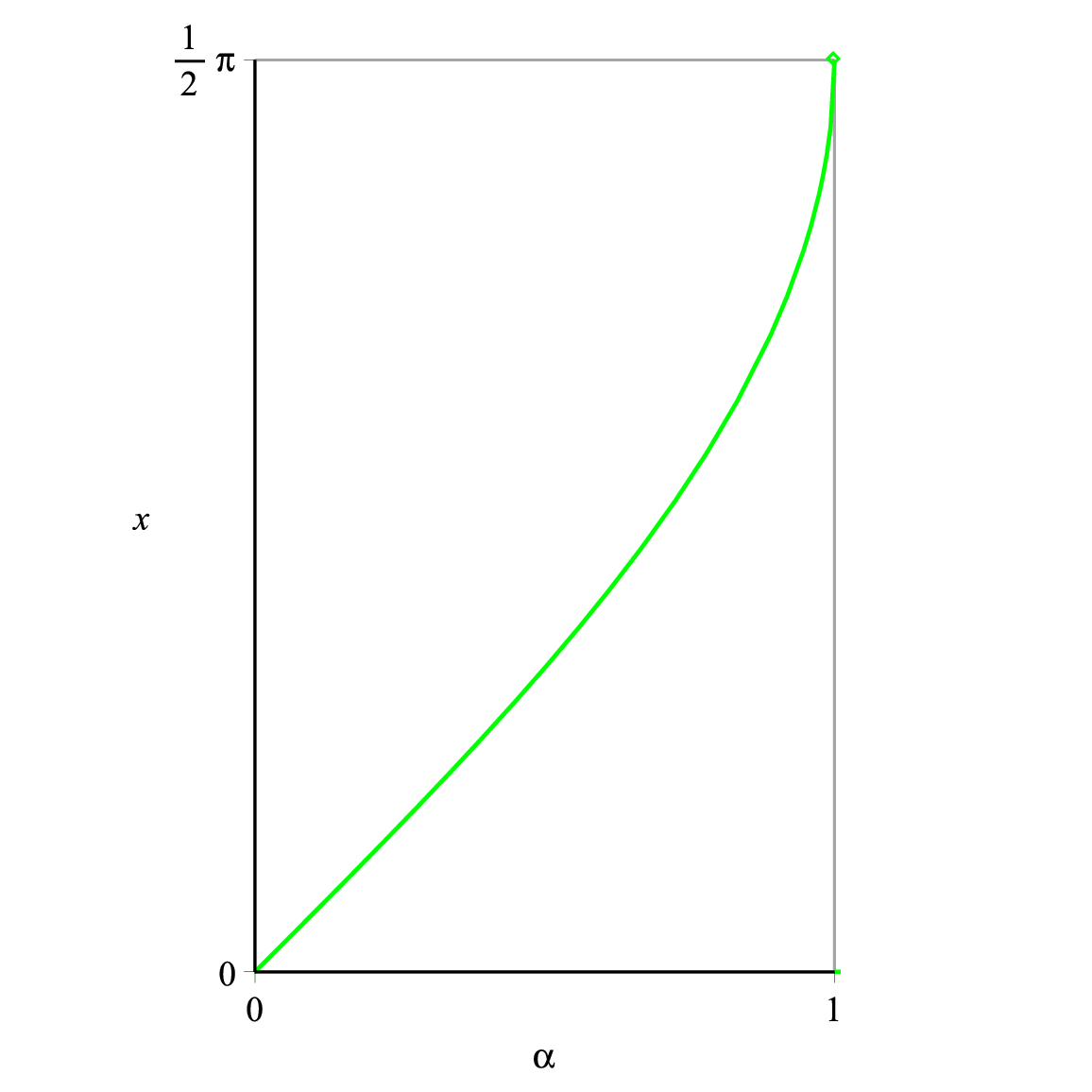}
\includegraphics[width=23mm]{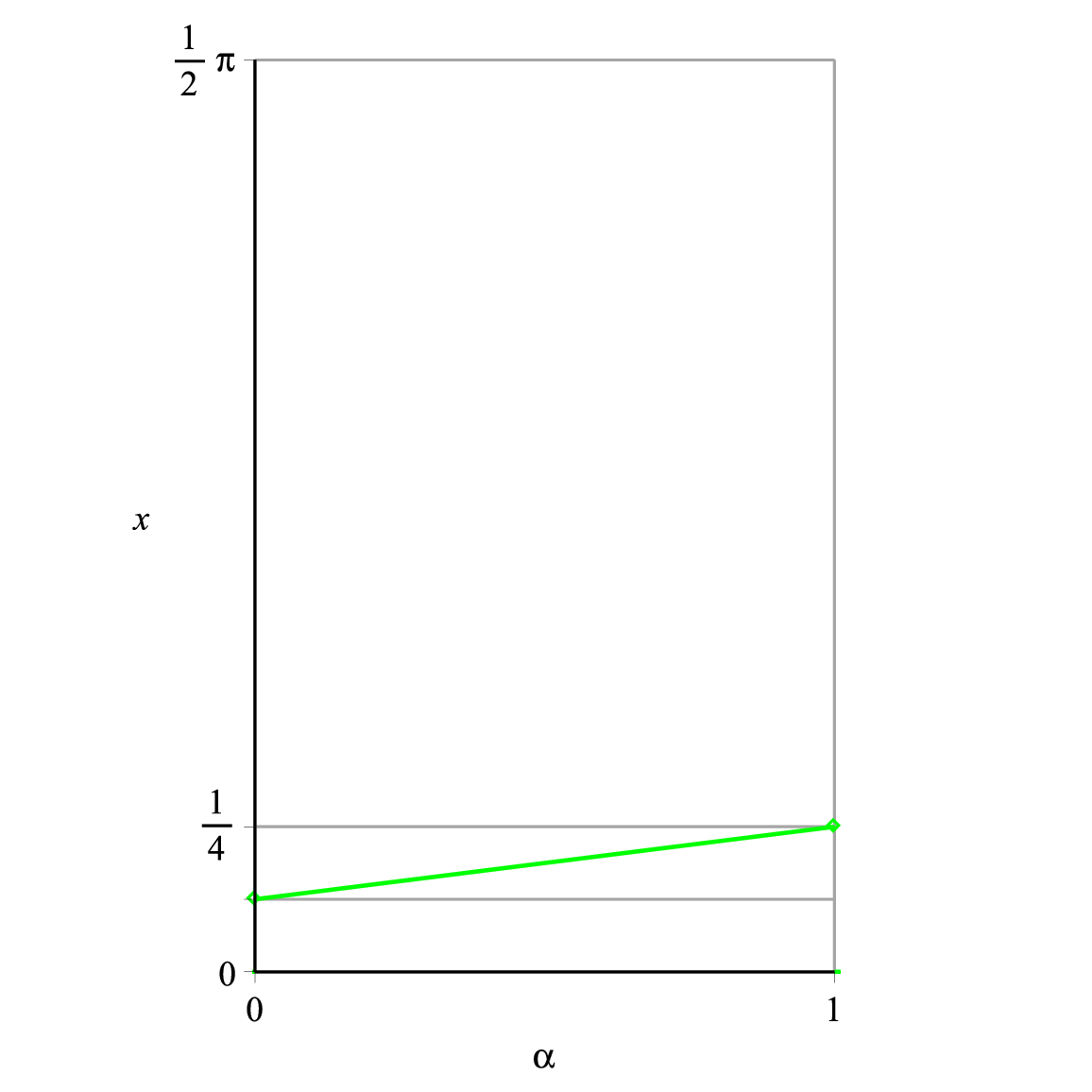}
\includegraphics[width=23mm]{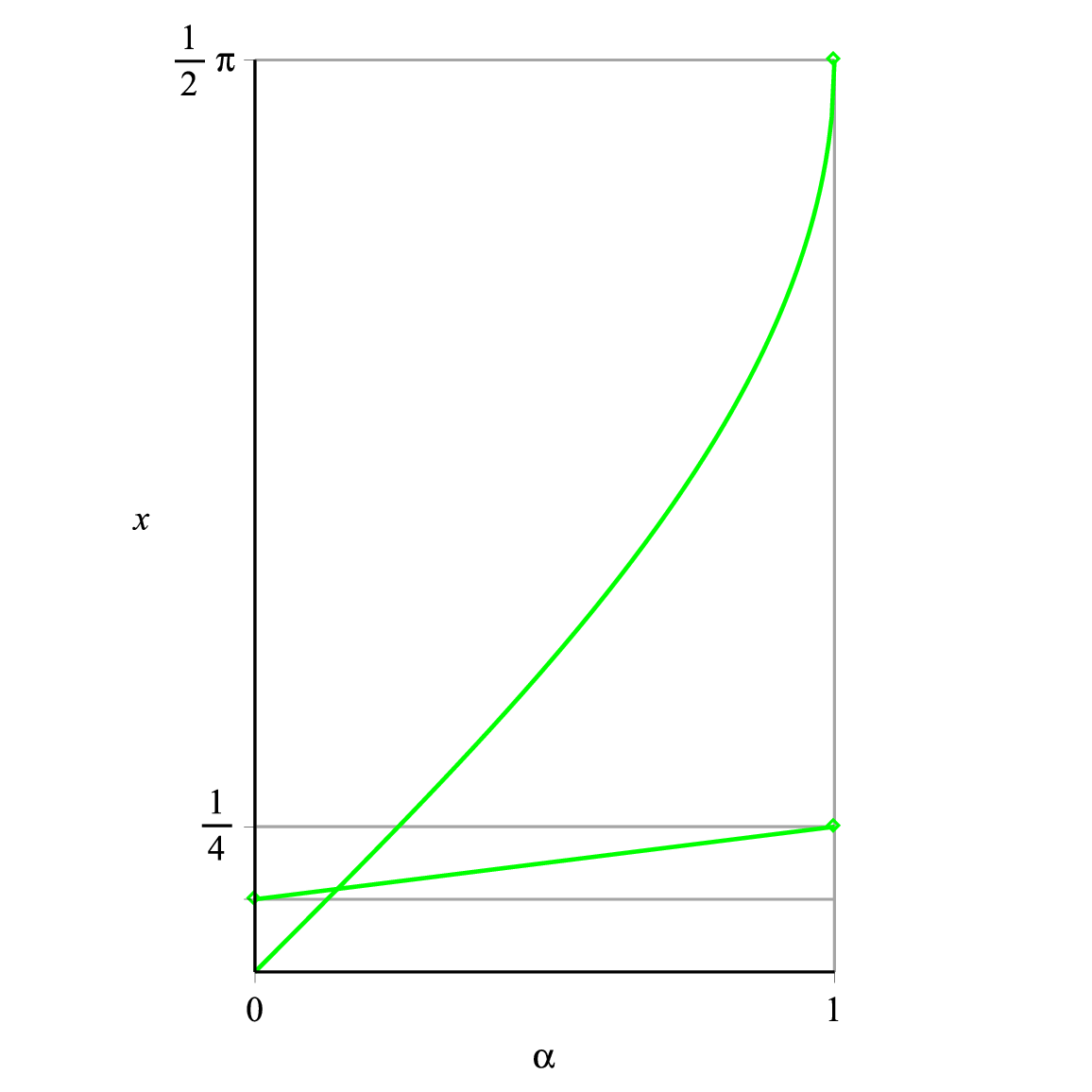}
\includegraphics[width=23mm]{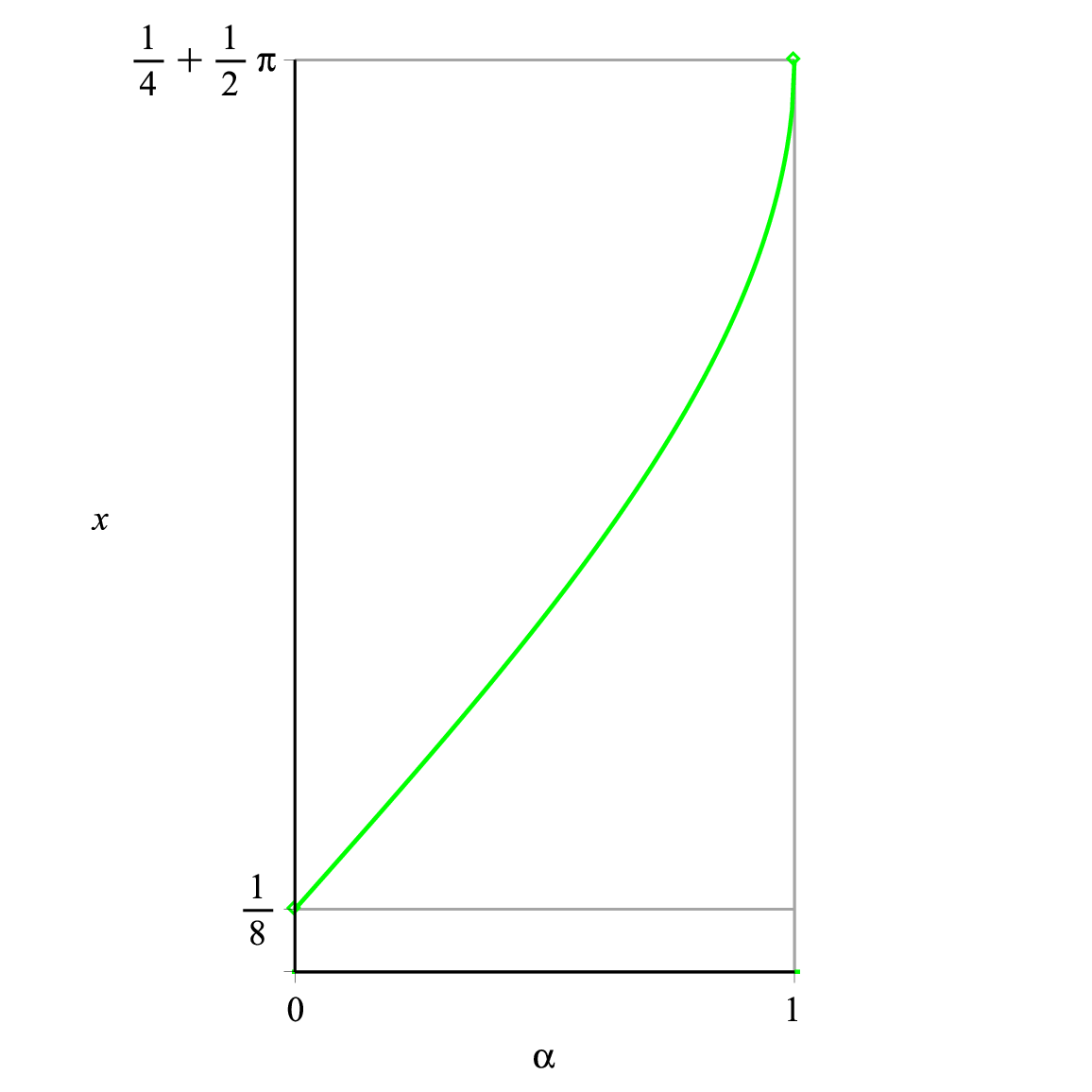}
\includegraphics[width=23mm]{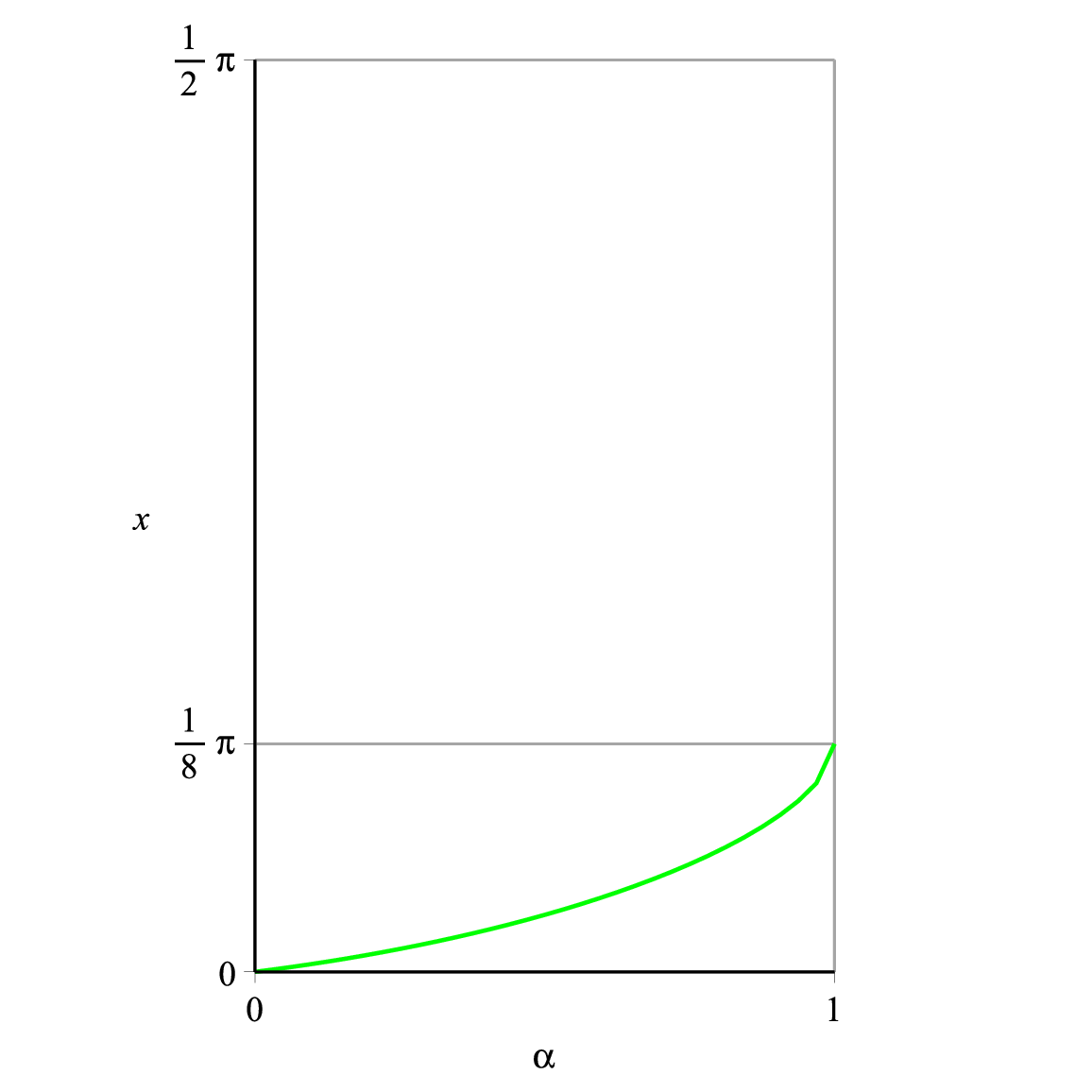}
\caption{\textit{From left to right:} \tiny $\xi_d(\alpha)$, $\eta_d(\alpha)$, $\xi_d(\alpha)$ and
$\eta_d(\alpha)$, $(\xi_d+\eta_d)(\alpha)$, $(\xi_d\cdot\eta_d)(\alpha)$.}
\end{center}
\end{figure}

\begin{figure}[H]
\begin{center}
\includegraphics[height=55mm]{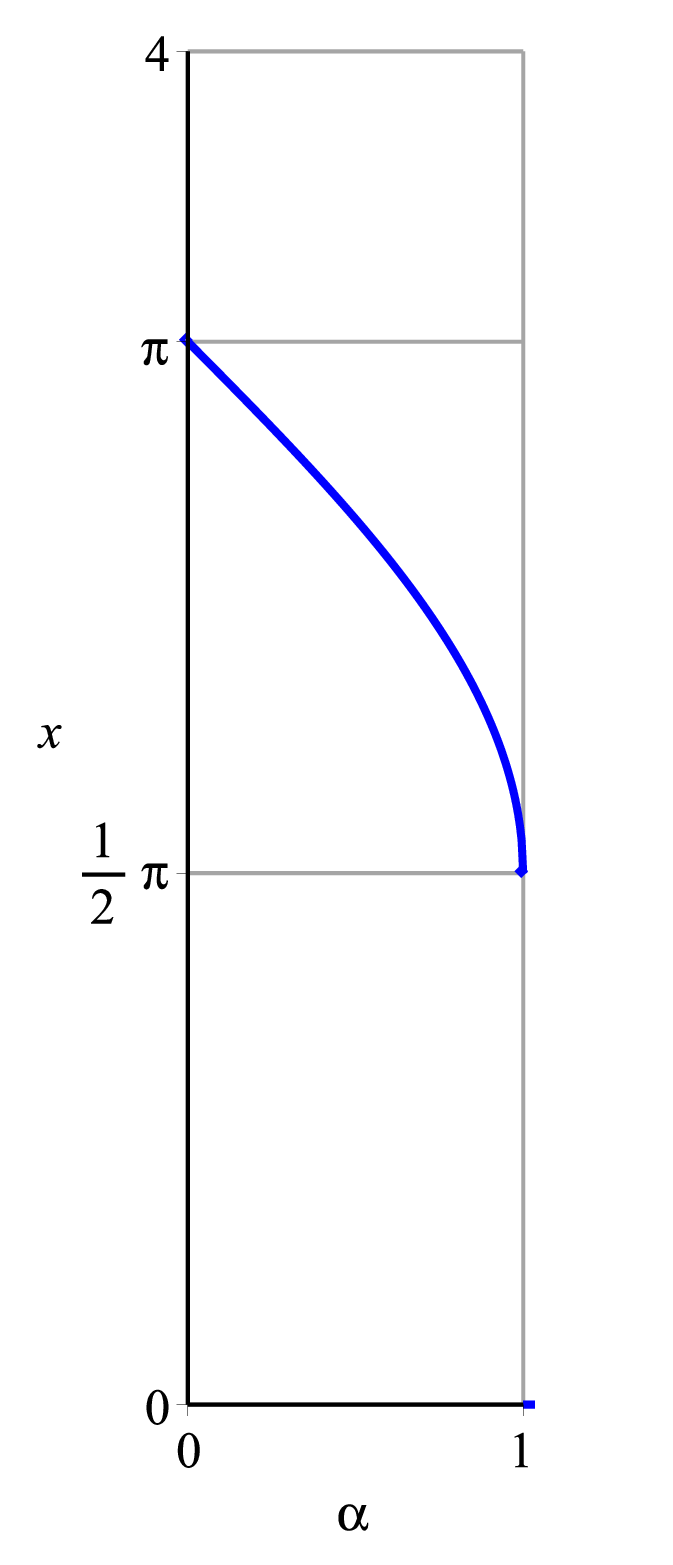}
\includegraphics[height=55mm]{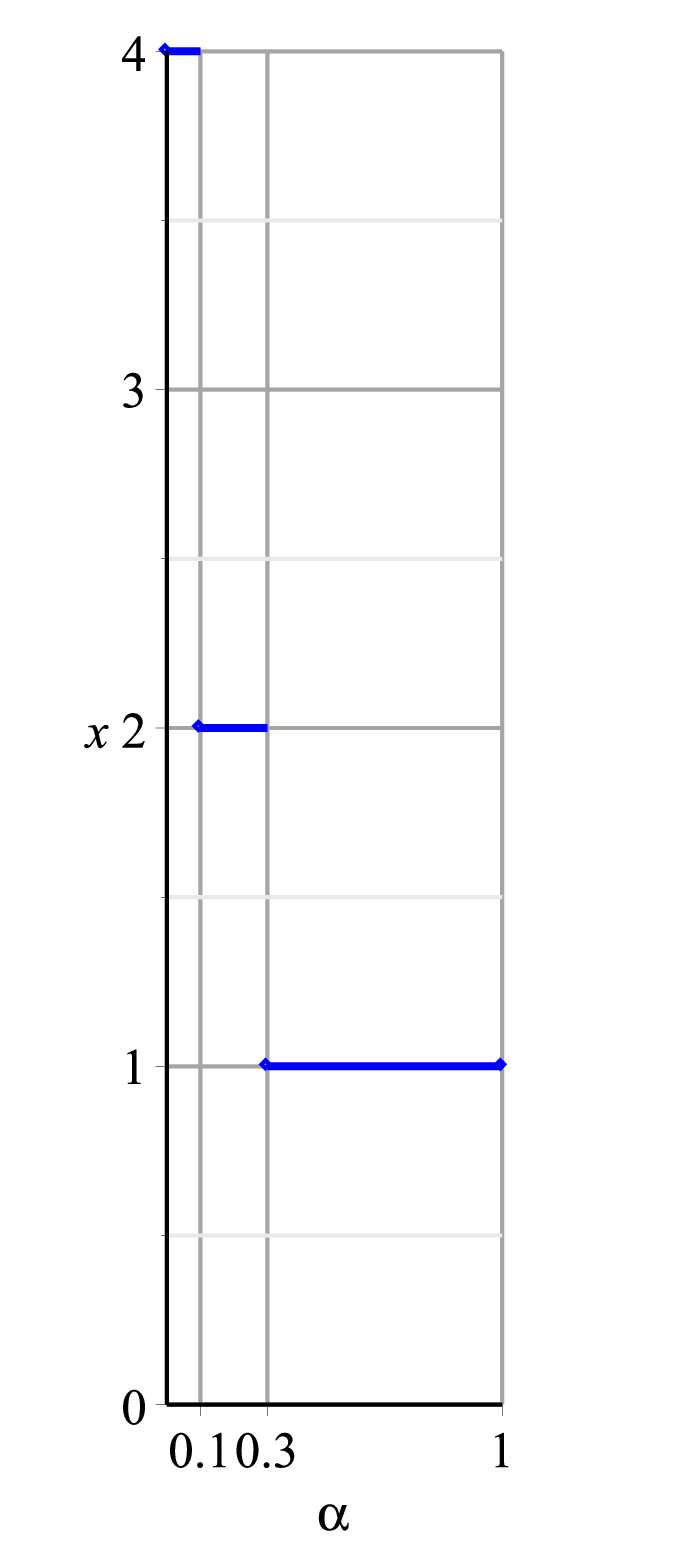}
\includegraphics[height=55mm]{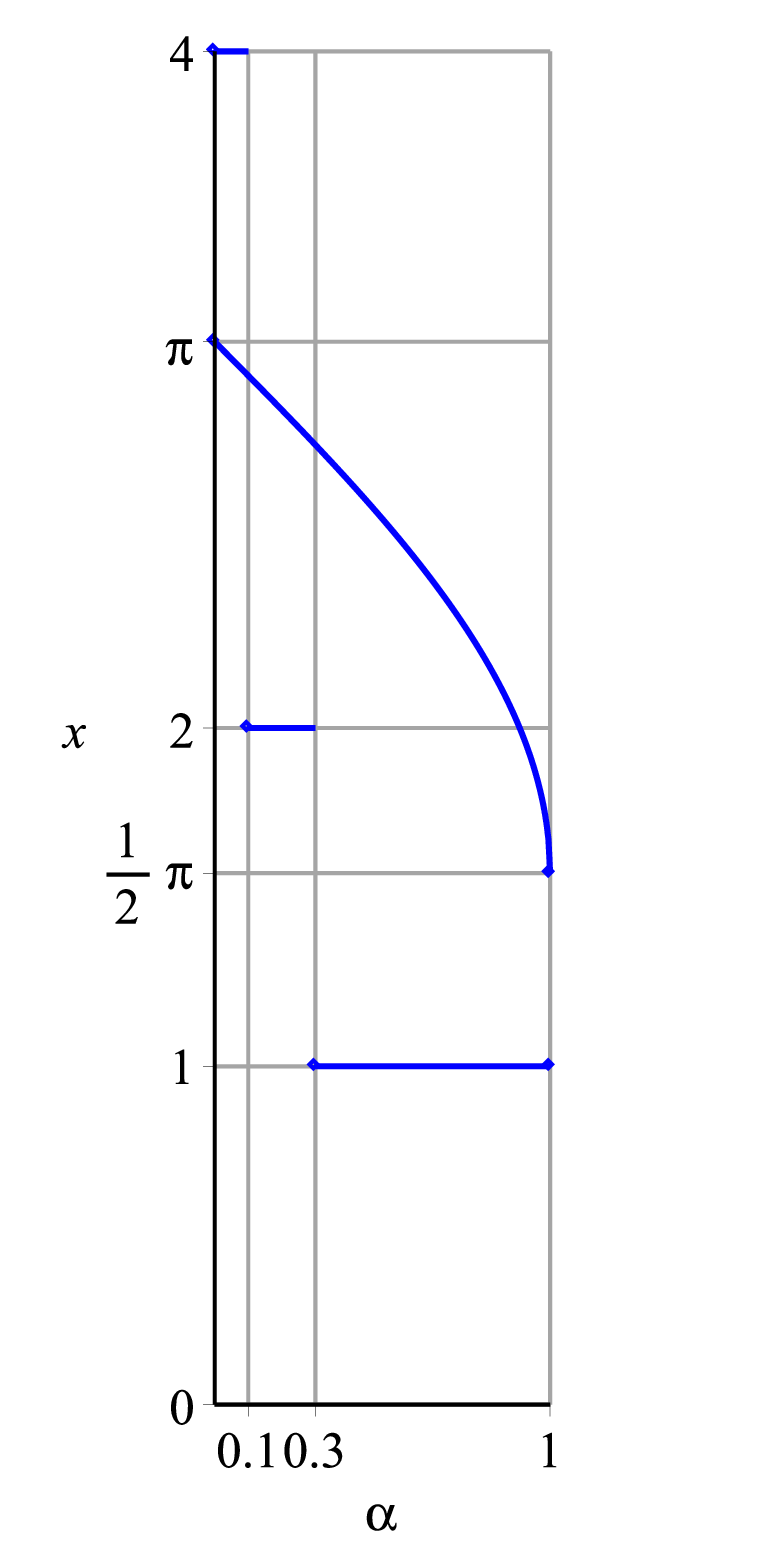}
\includegraphics[height=55mm]{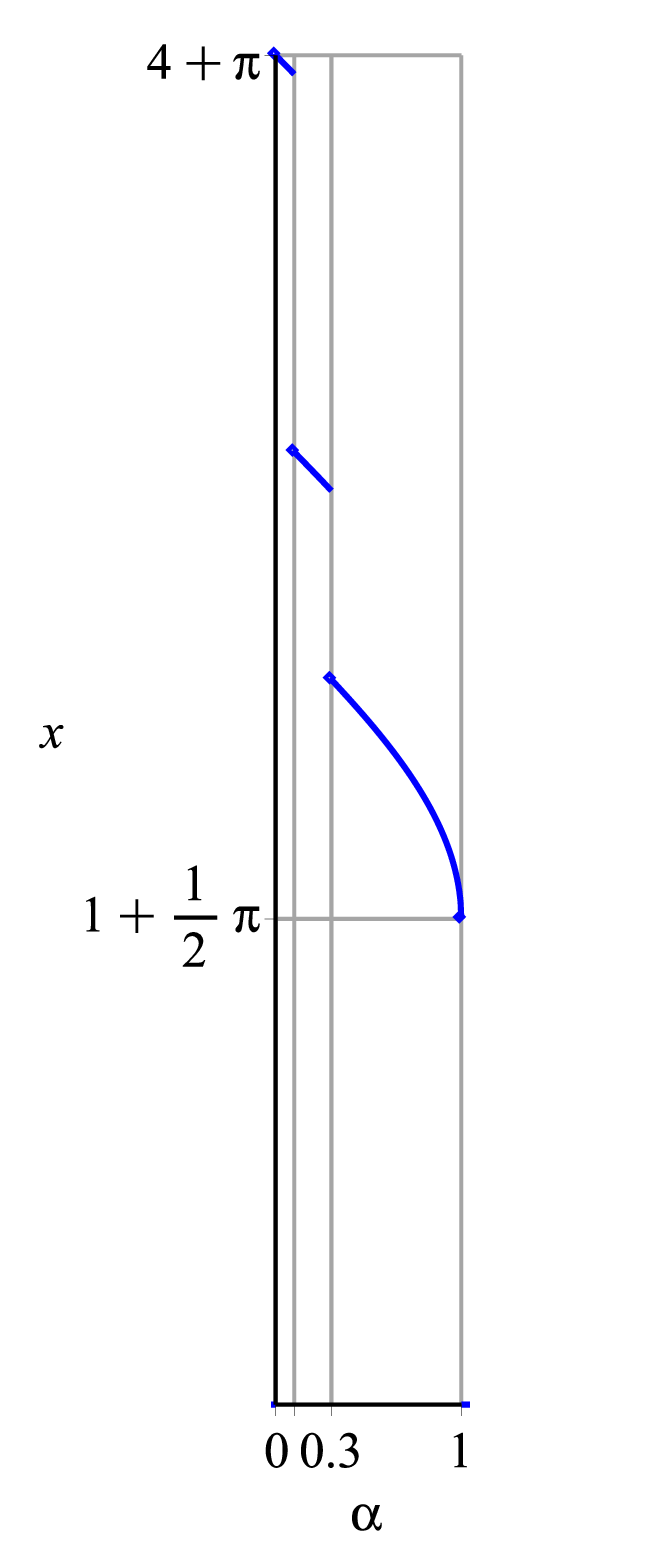}
\includegraphics[height=55mm]{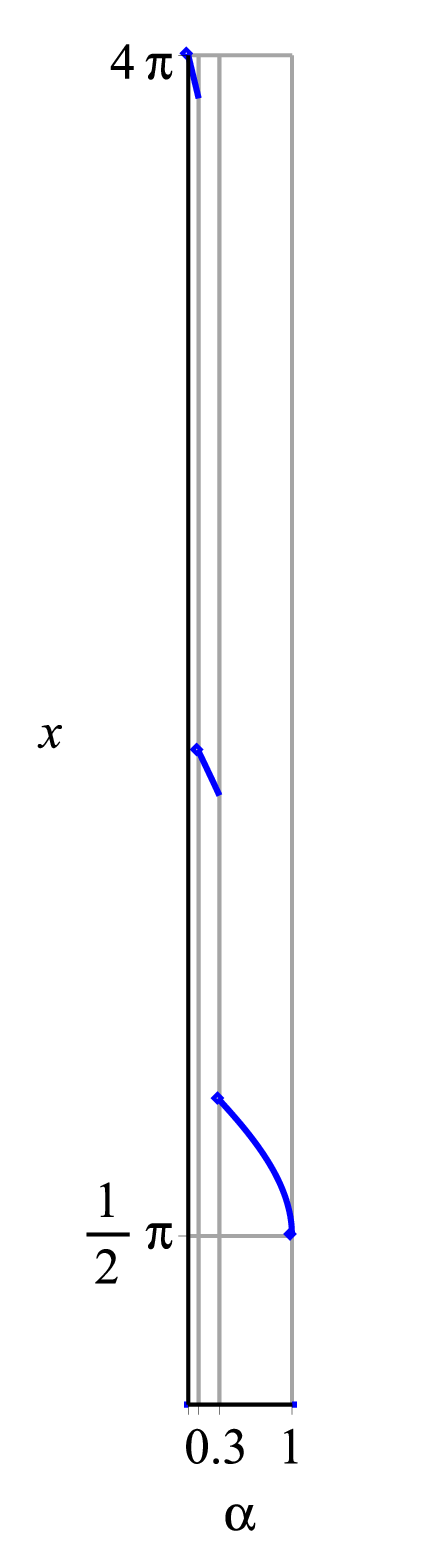}
\caption{\textit{From left to right:} \tiny $\xi_u(\alpha)$, $\eta_u(\alpha)$, $\xi_u(\alpha)$ and
$\eta_u(\alpha)$, $(\xi_u+\eta_u)(\alpha)$, $(\xi_u\cdot
\eta_u)(\alpha)$}
\end{center}
\end{figure}

After inverting back once more we obtain

\begin{equation*}
   \xi\odot\eta=
      \begin{cases}
      (\xi\odot\eta)_l=&  \{\arcsin(\alpha)\cdot(\frac{\alpha-1}{8}) = x\},\quad x\in[0,\frac{\pi}{8}],\,
      \alpha\in[0,1]\\
      (\xi\odot\eta)_m=&  [\frac{\pi}{8},\frac{\pi}{2}]\\
      (\xi\odot\eta)_r=&  \phi(x)\,,
      \end{cases}
\end{equation*}

where

\begin{equation*}
   \phi(x)=
     \begin{cases}
      \sin(1^{-1}x)\quad&x\in[(\pi-\arcsin(1))\cdot1,(\pi-\arcsin(0.3))\cdot1]\\
      0.3  &x\in((\pi-\arcsin(0.3)\cdot1,(\pi-\arcsin(0.3))\cdot2]\\
      \sin(2^{-1}x)&x\in((\pi-\arcsin(0.3))\cdot2,(\pi-\arcsin(0.1))\cdot2]\\
      0.1  &x\in ((\pi-\arcsin(0.1))\cdot2,(\pi-\arcsin(0.1))\cdot4]\\
      \sin(4^{-1}x)&x\in((\pi-\arcsin(0.1))\cdot4,(\pi-\arcsin(0))\cdot4]\,.
     \end{cases}
\end{equation*}\bigskip

Similarly

\begin{equation*}
   \xi\oplus\eta=
      \begin{cases}
      (\xi\oplus\eta)_l=&  \{\arcsin(\alpha)+(\frac{\alpha+1}{8}) = x\},\quad x\in[0,\frac{\pi}{2}+\frac18],\alpha\in[0,1]\\
      (\xi\oplus\eta)_m=&  [\frac{\pi}{2}+\frac18,\frac{\pi}{2}+1]\\
      (\xi\oplus\eta)_r=&  \phi(x)\,,
      \end{cases}
\end{equation*}

with

\begin{equation*}
   \phi(x)=
     \begin{cases}
      \sin(x-1)\quad &x\in[(\arcsin(1)+1,(\pi-\arcsin(0.3))+1]\\
      0.3 &x\in((\pi-\arcsin(0.3)+1,(\pi-\arcsin(0.3))+2]\\
      \sin(x-2) &x\in((\pi-\arcsin(0.3))+2,(\pi-\arcsin(0.1))+2]\\
      0.1 &x\in ((\pi-\arcsin(0.1))+2,(\pi-\arcsin(0.1))+4]\\
      \sin(x-4) &x\in((\pi-\arcsin(0.1))+4,(\pi-\arcsin(0))+4]\,.
     \end{cases}
\end{equation*}

\bigskip

\bigskip

Concluding the calculations of the particular components of the outcome let us notice that:\par
\begin{itemize}
\item
The components $(\xi\odot\eta)_l$ and $(\xi\oplus\eta)_l$ are in a form often encountered in practice, where the function is given in
implicit terms because providing a straightforward closed form representation is impossible.
\item
The functions $(\xi\odot\eta)_r$ and $(\xi\oplus\eta)_r$ are given in explicit terms where one could have applied~\eqref{dodawanie} directly or simpler  \eqref{mnozenie przez rzeczywista} and \eqref{dodawanie rzeczywistej} as always taking extreme care to ascertain the proper (open or closed) endpoints of intervals.
\item
The middle part can be always easily written out using \eqref{przedzialowe}.
\end{itemize}

\bigskip
Here are the graphs of $\xi$ and $\eta$ again:\bigskip

Before conducting the arithmetic operations

\begin{figure}[h]
\includegraphics[height=10mm]{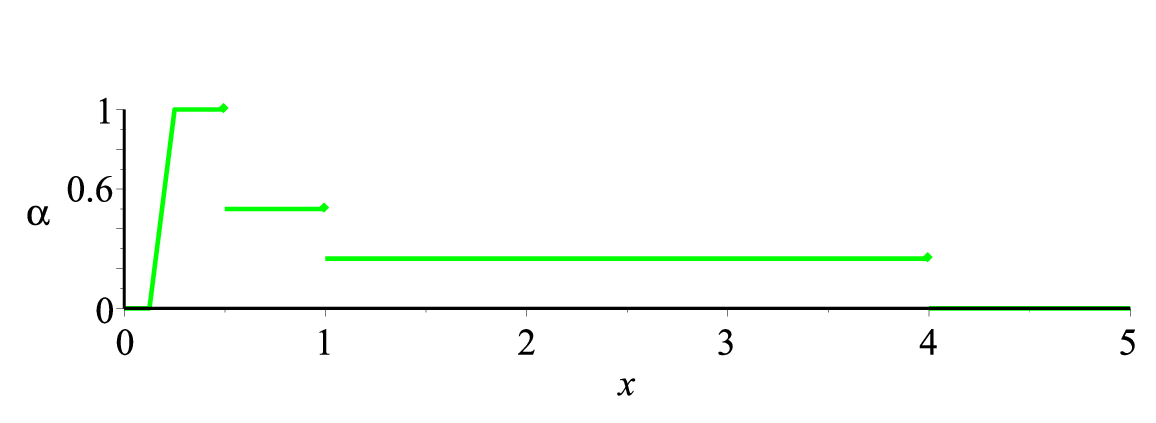}
\end{figure}
\begin{figure}[h]
\begin{center}
\includegraphics[height=10mm]{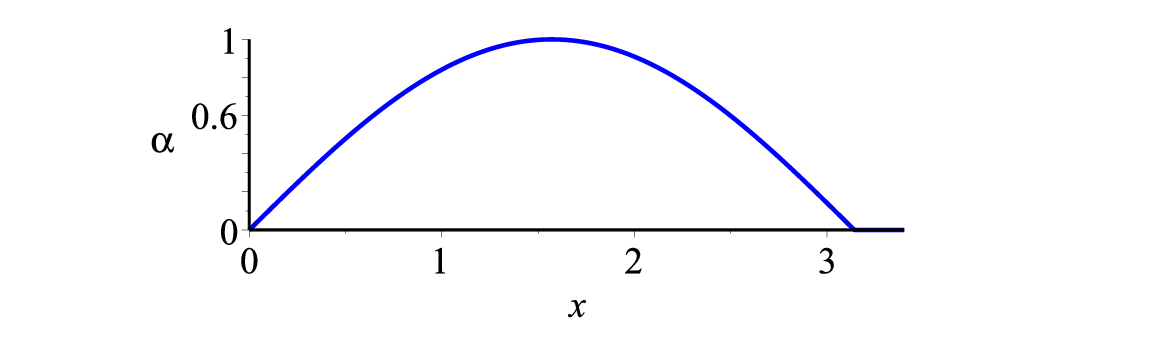}
\end{center}
\end{figure}

after addition
\begin{center}
\includegraphics[height=19mm]{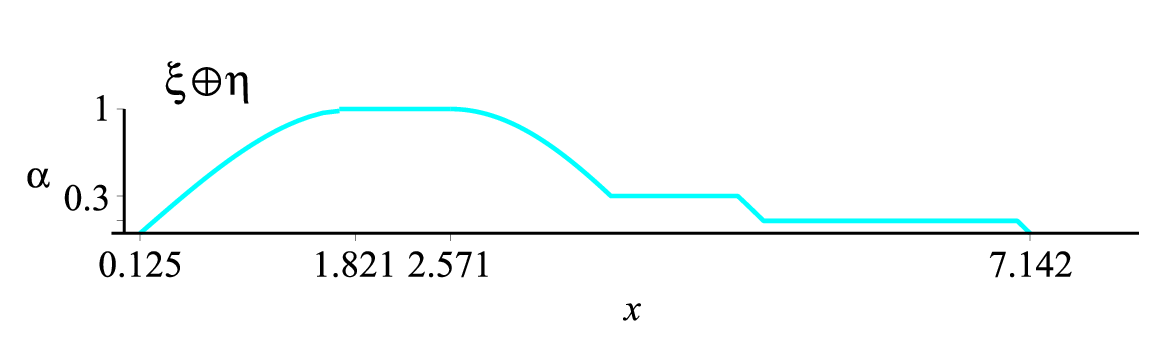}
\end{center}

and after multiplication
\begin{center}
\includegraphics[height=25mm]{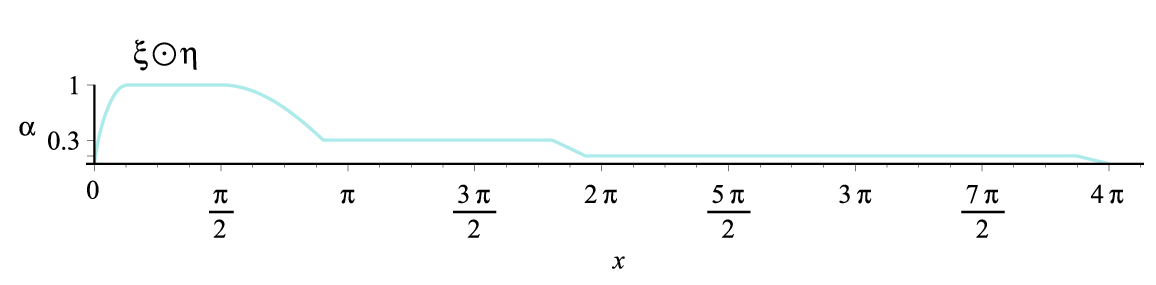}\,.
\end{center}
\end{example}

\begin{subsection}{Algebraic properties of the class $\Ft$}
\bigskip
From the general properties of functions and their inverses we see that the operations ``$\oplus$'' and ``$\odot$'' defined  on
$\Ft$ are commutative and associative.  We also have the zero element $\delta_0$ and the multiplicative identity
$\delta_1$, so $\Ft\cup\{\delta_0\}$ is a commutative semiring.

\noindent In general, the elements of $\Ft$ do not have inverse elements, which can be easily seen in the case of interval
numbers:\par

\smallskip
Indeed, for addition, using \eqref{przedzialowe} we would have to have
$$[a,b]\oplus[x,y]= [0,0],$$
so $x=-a, y=-b$. But $[-a,-b]$ with $a<b$ does not have the right orientation meaning that the called for interval $[x,y]$ does not exist (even when allowing for negatives, as we will in section~3).\par
\smallskip
For multiplication, using \eqref{przedzialowe} we would have to have
$$[a,b]\odot[x,y]= [1,1],$$
but $x=\frac{1}{a}\ge y=\frac{1}{b}$, meaning again there is no inverse in the strict, non-fuzzy sense. Note however that $[\frac{a}{b},\frac{b}{a}]$ will pass for a fuzzy $1$.
\end{subsection}

\newpage

\section{$\Fth$}
Throughout sections $1$ and $2$ we have been assuming the support of the characterizing functions to be a finite closed interval lying within $\R^+$.
This assumption greatly facilitates the understanding of the essence of how fuzzy arithmetic works. In \textit{Definition~3} below this assumption is dropped. Negative
and
mixed (i.e. containing zero) supports are also considered. \textit{Definition~3} is what one would normally encounter in a textbook and this particular
version
is taken from \cite{Viertl11}.
Subsequently the definition is translated into the language used in the preceding sections.

\bigskip
Here is the compact form:
\begin{definition}\label{def3a}
A fuzzy interval $\xi$ is determined by its characterizing function $\xi(\cdot)$
which is a real function of one real variable $x$ obeying the following
\begin{itemize}
\item[(0)] $\xi(\cdot):\R\mapsto[0,1].$
\item[(1)] $\exists x_0\in\R$ such that $\xi(x_0)=1.$
\item[(2)] $\xi(\lambda x_1+(1-\lambda)x_2)\geq\min\left(\xi(x_1),\xi(x_2)\right).$ (fuzzy-convexity)
\item[(3)] $\xi(\cdot)$ is upper semi-continuous.\quad($\lim_{x_n\to x_0}\xi(x_n)\leq\xi(x_0)$)
\item[(4)] $\xi(\cdot)$ has compact support.
\end{itemize}
\end{definition}
\smallskip
There is an equivalent (for a proof of this see \cite{Viertl11}) and very convenient definition formulated in the widespread language of so-called
\textit{$\alpha$-cuts} which a physicist would rather refer to as level-sets or isolines:\par\vspace{1mm}
\noindent\textbf{Definition 3'.}\label{def3'}
A fuzzy interval $\xi$ is determined by its characterizing function $\xi(\cdot)$ which is a real function of one real variable $x$ obeying the following:
\begin{itemize}
\item[(a)] $\xi : \R \mapsto [0, 1].$
\item[(b)] $\forall \alpha\in(0,1]$ the so-called $\alpha$-cut $C_\alpha(\xi)=\{x\in\R:\xi(x)\geq\alpha\}$ is a (non-empty) compact interval.
\item[(c)] The support of $\xi(\cdot)$, $\supp[\xi(\cdot)]:= \overline{\{x\in\R : \xi(x) > 0\}}$ is bounded.
\end{itemize}

\begin{remark}
Note that $C_\alpha(\xi)=\{x\in\R:\xi(x)\geq\alpha\}=[{\xi_d}(\alpha),{\xi_u}(\alpha)]$, and this is the bridge to how we have been operating up to this moment.
\end{remark}

\textbf{Definition 3''}
As in sections $1$ and $2$ we may equivalently define a fuzzy interval (number) $\xi$ to be defined by an ordered pair $\left(\xi_l(\cdot),\xi_r(\cdot)\right)$ of two real functions
$\xi_l(\cdot):\R\to[0,1]$ and $\xi_r(\cdot):\R\to[0,1]$ of one real variable $x$ such that:\par\noindent
\begin{itemize}
\item[(1)] $\xi_l(\cdot)$ is increasing, right-continuous and has support on some closed interval $[l,\underline{m}]$
\item[(2)] $\xi_r(\cdot)$ is decreasing, left-continuous and has support on some closed interval $[\overline{m},r]$ such that
    $\underline{m}\leq\overline{m}$
\item[(3)] $\xi_l(\underline{m})=\xi_r(\overline{m})=1$
\item[(2')] To differentiate fuzzy numbers from fuzzy intervals we additionally request: $\xi_d(1)=\{\underline{m}\}$, and $\xi_u=\{\overline{m}\}$.
\end{itemize}

\bigskip
Just like in real (not fuzzy) interval arithmetic we define for all four arithmetic operations:

\begin{equation}
{(\xi\circ\eta)}_d:= \min(\xi_d\circ\eta_d,\xi_d\circ\eta_u,\xi_u\circ\eta_d,\xi_u\circ\eta_u)
\end{equation}
and
\begin{equation}
{(\xi\circ\eta)}_u:= \max(\xi_d\circ\eta_d,\xi_d\circ\eta_u,\xi_u\circ\eta_d,\xi_u\circ\eta_u),
\end{equation}

\smallskip
with ``$\circ$'' standing for any of the four operations: ``$\oplus$'',``$\ominus$'',``$\odot$'',``$\oslash$''.\par
\bigskip
Then
\begin{equation}
(\xi\circ\eta)_l:=\overleftarrow{{(\xi\circ\eta)_d}}\,,
\end{equation}
and
\begin{equation}
(\xi\circ\eta)_r:=\overleftarrow{(\xi\circ\eta)_u}\,.
\end{equation}
\smallskip

And similar to (\ref{lrform2}) before

\begin{equation}
\xi\circ\eta=
\begin{cases}
\overleftarrow{{(\xi\circ\eta)_d}}\,,\\
1\\
\overleftarrow{{(\xi\circ\eta)_u}}\,.
\end{cases}
\end{equation}

\bigskip
Or more compact, in interval notation:
\begin{equation}\label{generaloperation}
\xi\circ\eta=\left((\xi\circ\eta)_l,(\xi\circ\eta)_r\right)=\left(\overleftarrow{{(\xi\circ\eta)_d}}\,,\overleftarrow{{(\xi\circ\eta)_u}}\right)
\end{equation}

\bigskip
We now separately turn our attention to each of the four arithmetic operations:

\begin{subsection}{The four operations}
\begin{subsubsection}{Addition}\label{addition3}
Because real addition is monotone in both variables we find that $(\xi\oplus\eta)_d=\xi_d+\eta_d$ and $(\xi\oplus\eta)_u=\xi_u+\eta_u$ and so
\eqref{generaloperation} reduces to \eqref{uododawanie}:
\begin{equation}
\xi\oplus\eta=\left((\xi\oplus\eta)_l,(\xi\oplus\eta)_r\right)=\left(\overleftarrow{(\xi_d+\eta_d)},\overleftarrow{(\xi_u+\eta_u)}\right)
\end{equation}

\end{subsubsection}

\begin{subsubsection}{Subtraction}
Define $$\xi\ominus\eta:=\xi\oplus-\eta$$ where $-\eta=-(\eta_l,\eta_r)=(-\eta_r,-\eta_l)$\par\noindent

\smallskip
Then fuzzy subtraction is defined as addition above.
\begin{remark}
Note that $\xi\ominus\xi=(\xi_l-\xi_r,\xi_r-\xi_l)$ is a symmetrical fuzzy zero but certainly not $\delta_0$. (Its support is
$\supp[\xi\ominus\xi]=[l-r,r-l]$).
\end{remark}
 \end{subsubsection}

\begin{subsubsection}{Mupltiplication}
\begin{equation}\label{multiplication3}
\begin{split}
& \xi\odot\eta =\\ &=\biggl(\overleftarrow{(\min(\xi_d\cdot\eta_d,\xi_d\cdot\eta_u,\xi_u\cdot\eta_d,\xi_u\cdot\eta_u)}\:,\overleftarrow{(\max(\xi_d\cdot\eta_d,\xi_d\cdot\eta_u,
\xi_u\cdot\eta_d,\xi_u\cdot\eta_u)}\biggr)
\end{split}
\end{equation}
Unfortunately, because unlike in the case of addition the real multiplication operator is not monotone, formula \eqref{generaloperation} does not reduce
automatically to \eqref{uomnozenie} and calculations can become very laborious and error prone when done by hand.
At the very end of section 4 some examples of multiplied triangle numbers of mixed support are graphically illustrated.\par
\smallskip
\begin{remark}\label{multiplication3'}
If $\xi,\eta\in\Ft$ \eqref{uomnozenie} applies as demonstrated in section 2.
\end{remark}
\end{subsubsection}

\begin{subsubsection}{Division}
Division is defined as the fuzzy inverse of multiplication by setting
\begin{equation}\label{division3}
[(\xi\oslash\eta)_d,(\xi\oslash\eta)_u)]=\biggl[\min(\frac{\xi_d}{\eta_d},\frac{\xi_d}{\eta_u},\frac{\xi_u}{\eta_d},\frac{\xi_u}{\eta_u}),\max(\frac{\xi_d}{\eta_d},\frac{\xi_d}{\eta_u},
\frac{\xi_u}{\eta_d},\frac{\xi_u}{\eta_u})\biggr]\,.
\end{equation}
For the above to be well defined we must assume $0\notin\supp{\eta_l}\cup\supp{\eta_r}.$
\begin{remark}
For $\eta_d$ left-continuous $1\slash{\eta_d}$ is right-continuous and likewise for $\eta_u$ right-continuous we have $1\slash{\eta_u}$ left continuous, so
the quotients in \eqref{division3} are all semi-continuous ``in the right way".\par\smallskip
\end{remark}
\begin{remark}\label{positivedivision}
For $\xi,\eta\in\Ft$ \eqref{division3} becomes just $$[(\xi\oslash\eta)_d,(\xi\oslash\eta)_u)]=\left[\frac{\xi_d}{\eta_u},\frac{\xi_u}{\eta_d}\right]\,.$$
\end{remark}
\bigskip
Let us now return to our very first worked through example 1 from section~1, $\xi$ as in \eqref{xiexample1}, $\eta$ as in \eqref{etaexample1}:
\begin{example}{(\textbf{fuzzy division})}
Going by \eqref{division3} and \textit{Remark} \ref{positivedivision}:
\begin{equation}
\frac{\xi_d}{\eta_u}=\frac{{\xi_l}^{-1}}{{\eta_r}^{-1}}=\frac{\alpha+1}{10-3\alpha} \text{ for }\alpha\in[0,1]\,,
\end{equation}
and hence
\begin{equation}
{(\xi\oslash\eta)}_l=\frac{10x-1}{1+3x}\text{ for }x\in[1,2]\oslash[7,10]=\left[\frac{1}{10},\frac{2}{7}\right]\,.
\end{equation}
Likewise
\begin{equation}
\frac{\xi_u}{\eta_d}=\frac{{\xi_r}^{-1}}{{\eta_l}^{-1}}=\frac{3-\alpha}{2\alpha+5} \text{ for }\alpha\in[0,1]\,,
\end{equation}
and hence
\begin{equation}
{(\xi\oslash\eta)}_r=\frac{3-5x}{2x+1}\text{ for }x\in[2,3]\oslash[5,7]=\left[\frac27,\frac35\right]\,.
\end{equation}
\begin{figure}[H]
\includegraphics[width=120mm]{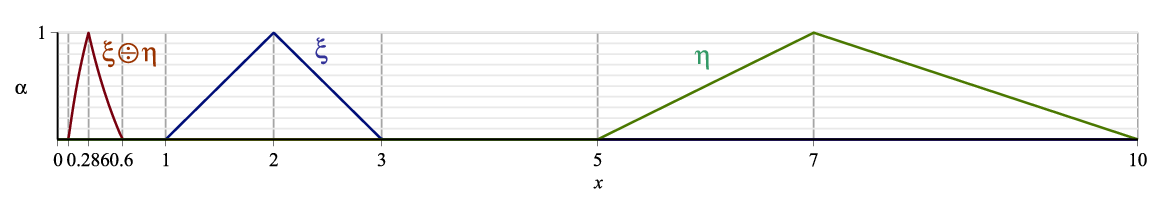}
\caption{Fuzzy division of two triangle numbers}
\end{figure}
\end{example}
\end{subsubsection}
\end{subsection}

\bigskip
\begin{section}{$\Fth^{-1}$ - reversing the point of view}

Especially in the case when fuzzy intervals are actively designed to model an economic situation to begin with as in management science, as opposed to being an expression of uncertainty in measurement of the outcome of an experiment in a natural science, it is obviously more natural and certainly more straightforward for computation to start with the right-hand sides of equations
\eqref{geninvinf} and \eqref{geninvsup} as primary from the very outset.

This leads to the so-called parametric representation of fuzzy intervals as suggested by R. Goetschel and W. Voxmann in~\cite{GV83} and~\cite{GV86}, and which has been predominant in the recent literature.

Therefore we propose to adopt the following definition, and associated terminology:\par\noindent
\begin{definition}
We understand a fuzzy (generalized) interval $$A^*\in[0,1]^\R\times[0,1]^\R$$ with \textit{fuzzy endpoints} $a_d(\cdot),a_u(\cdot)$ to be an ordered pair of functions in shorthand denoted by
$A^*:=[a_{d}(\cdot),a_{u}(\cdot)]$ such that:
\begin{itemize}
\item[(1)] $a_d(\cdot):[0,1]\mapsto\R$ is increasing and left-continuous,
\item[(2)] $a_u(\cdot):[0,1]\mapsto\R$ is decreasing, right-continuous.
\end{itemize}
Furthermore we demand in consistence with \textit{Definition 2(2)} and \textit{Remark~\ref{invsemicont}}:
\begin{itemize}
\item[(1')] $\lim_{\alpha\to{0^+}}a_d(\alpha)=a_d(0),$
\item[(2')] $\lim_{\alpha\to{1^-}}a_u(\alpha)=a_u(1).$
\end{itemize}

\end{definition}

\textit{Notation and Terminology:}
\begin{itemize}
\item
To distinguish real intervals $A$ from fuzzy intervals $A^*$ we use the starred form for the fuzzy interval.
\item
We as before call the (parameterizing) characterizing functions $a_d(\cdot)$, $a_u(\cdot)$ the \emph{fuzzy endpoints} of the fuzzy interval $[a_d,a_u]$ but use square brackets for the interval in parameterized form to distinguish it from the same interval in what we now might call the $x-$axis representation $\bigl(\overleftarrow{a_d}(\cdot),\overleftarrow{a_u}(\cdot)\bigr).$
\item
We propose to denote $a_d(1)=\underline{a}$ and $a_u(1)=\overline{a}$ and refer to the real interval $A=[\underline{a},\overline{a}]$ as the (real) \emph{core} (that is the real interval we are generalizing) of $A^*$.
\item
The \emph{support} of the fuzzy interval $A^*$ in parametric form is $\supp{A^*}=[d,u]$ with \mbox{$d=\alpha_{d}(0)$} and \mbox{$u=a_u(0)$.}
\end{itemize}

\begin{subsection}{Arithmetic operations on $\Fth^{-1}$}
For two fuzzy intervals $A^*=[a_d,a_u]$ and $B^*=[b_d,b_u]$ all four arithmetic operations on $\Fth^{-1}$ are then defined simply by:
\begin{equation}\label{generaloperation'}
\begin{split}
&A^*\circ B^*=\\
&\bigl[\min(a_d\circ b_d,a_d\circ b_u, a_u\circ b_d, a_u\circ b_u),\max(a_d\circ b_d, a_d\circ b_u,a_u\circ b_d,
a_u\circ b_u)\bigr]
\end{split}
\end{equation}
as in crisp interval arithmetic with ``$\circ$'' denoting any of the four operations: ``$\oplus,\ominus,\odot,\oslash$''.\par
\smallskip
All comments and simplifications of section 3 apply with the appropriate modifications.

\begin{example}
Real intervals $[\underline{a},\overline{a}]$ correspond to fuzzy endpoints\par
\vspace{1mm}$a_d(\alpha)=\underline{a} \text { for } \alpha\in[0,1]$ and\par
$a_u(\alpha)=\overline{a} \text { for } \alpha\in[0,1].$
\end{example}

\begin{example}
Fuzzy numbers are then fuzzy intervals in the sense above for the case the interval $[\underline{a},\overline{a}]$ is degenerated, that is
$\underline{a}=\overline{a}=x_0$ for some $x_{0}\in\R$.
\end{example}

\begin{example}
Triangle numbers $tr_{(l,m,r)}$ were characterized in section 1 by
\begin{equation}
tr(x)=
\begin{cases}
tr_l(x)=\frac{x-l}{m-l} &\text{ for } x\in[l,m]\\
1 &\text{ for } x=m\\
tr_r(x)=\frac{r-x}{r-m} &\text{ for } x\in[m,r]\,.
\end{cases}
\end{equation}
We now write simply
\begin{equation}
tr^*_{(d,x_0,u)}(\alpha)=[\alpha(x_0-d)+d,u-\alpha(u-x_0)], \quad\alpha\in[0,1]\,.
\end{equation}
Rewriting in parametric representation our very first Example~1 in which we considered the two triangle numbers \mbox{$tr_1=tr_{(1,2,3)}$} and \mbox{$tr_2=tr_{(5,7,10)}$}:
\begin{equation*}
\begin{split}
tr^*_1(\alpha)&=[\alpha(2-1)+1,3-\alpha(3-2)]=[\alpha+1,3-\alpha]\,,\\
tr^*_2(\alpha)&=[\alpha(7-5)+5,10-\alpha(10-7)]=[2\alpha+5,10-3\alpha].
\end{split}
\end{equation*}
\smallskip
Then by \eqref{generaloperation'}, \textit{Remark~\ref{multiplication3'}} and \textit{Remark~\ref{positivedivision}} we obtain
\begin{eqnarray}
\nonumber tr^*_1\oplus  tr^*_2 &=& \left[(\alpha+1)+(10-3\alpha),(3-\alpha)+(2\alpha+5)\right]=[10-2\alpha,8-3\alpha]\,,\\
\nonumber tr^*_1\ominus tr^*_2 &=& \left[(\alpha+1)-(2\alpha+5),(3-\alpha)-(10-3\alpha)\right]=[-\alpha-4,-7+2\alpha]\,,\\
\nonumber tr^*_1\odot   tr^*_2 &=& \left[(\alpha+1)\cdot(2\alpha+5),(3-\alpha)\cdot(10-3\alpha)\right]\,,\\
\nonumber tr^*_1\oslash tr^*_2 &=& \left[\frac{\alpha+1}{10-3\alpha},\frac{3-\alpha}{2\alpha+5}\right].
\end{eqnarray}
\noindent very simply and straightforward, without having to go through the hassle of inverting and re-inverting.
\end{example}
\end{subsection}

\smallskip
To recapitulate:\par
\noindent The main motivation and advantage of this section's approach is: By changing the point of view right from the start, when we set up a model and in doing so define a set of fuzzy intervals, the cumbersome procedure of (twice!) inverting normally un-invertible functions is
avoided.\par

This is especially true when characterizing functions of mixed support are involved, and we close this section with some random impressions of various multiplied triangle numbers of mixed support, ($tr^*_1$: green, $tr^*_2$: blue, $tr^*_1\odot tr^*_2$: red) that have been rewritten as in \textit{Example~15}:
\begin{figure}[H]
\includegraphics[width=41mm]{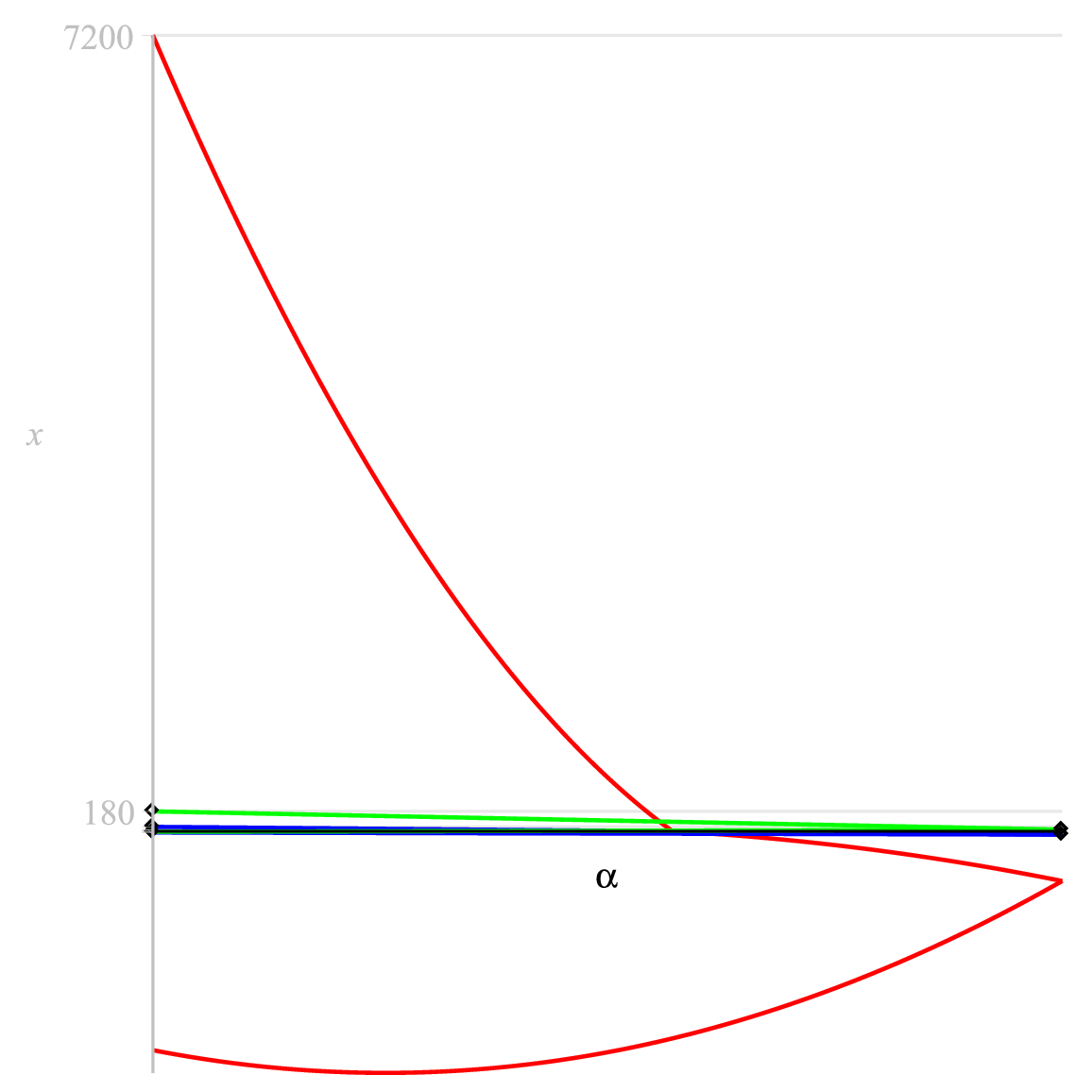}
\includegraphics[width=41mm]{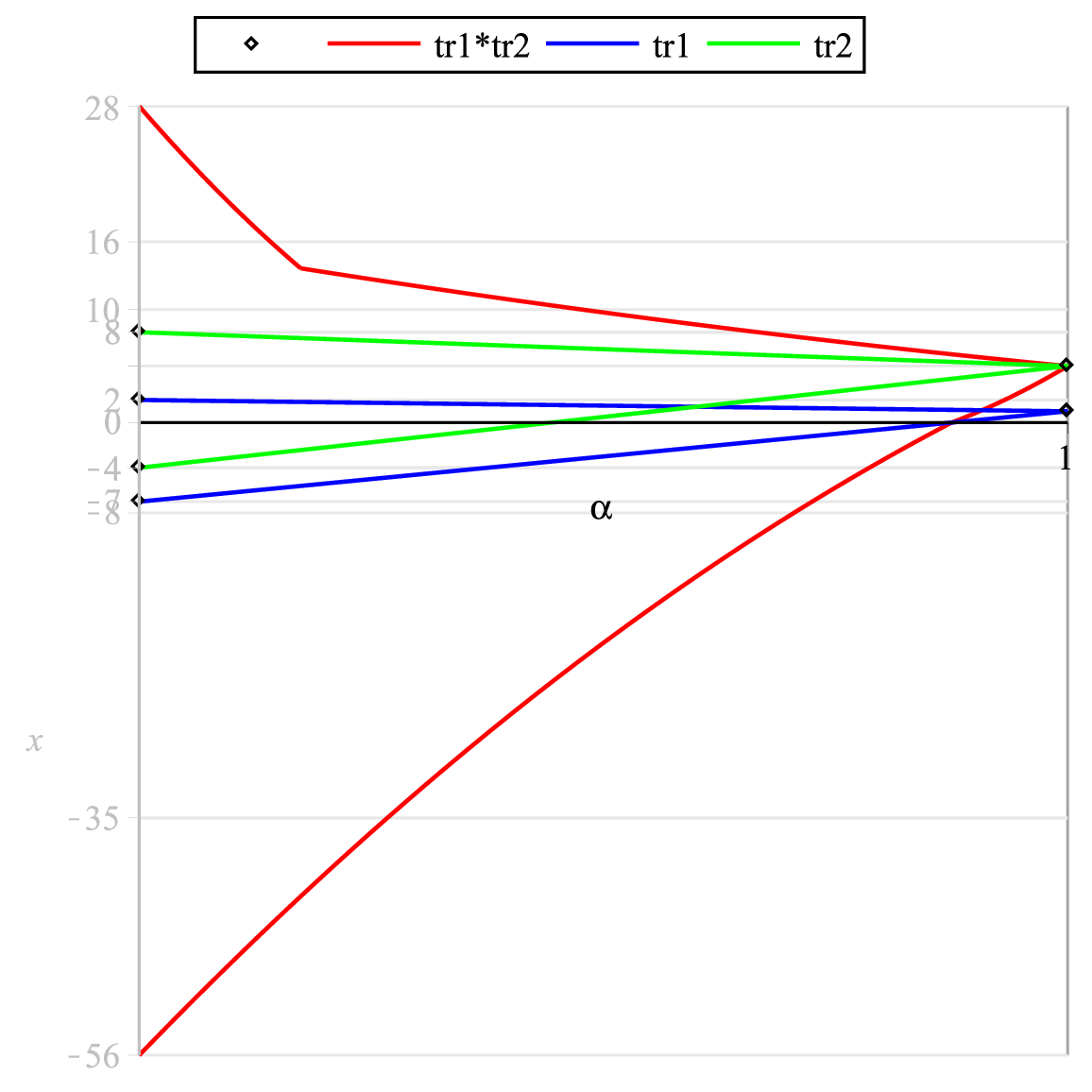}
\includegraphics[width=41mm]{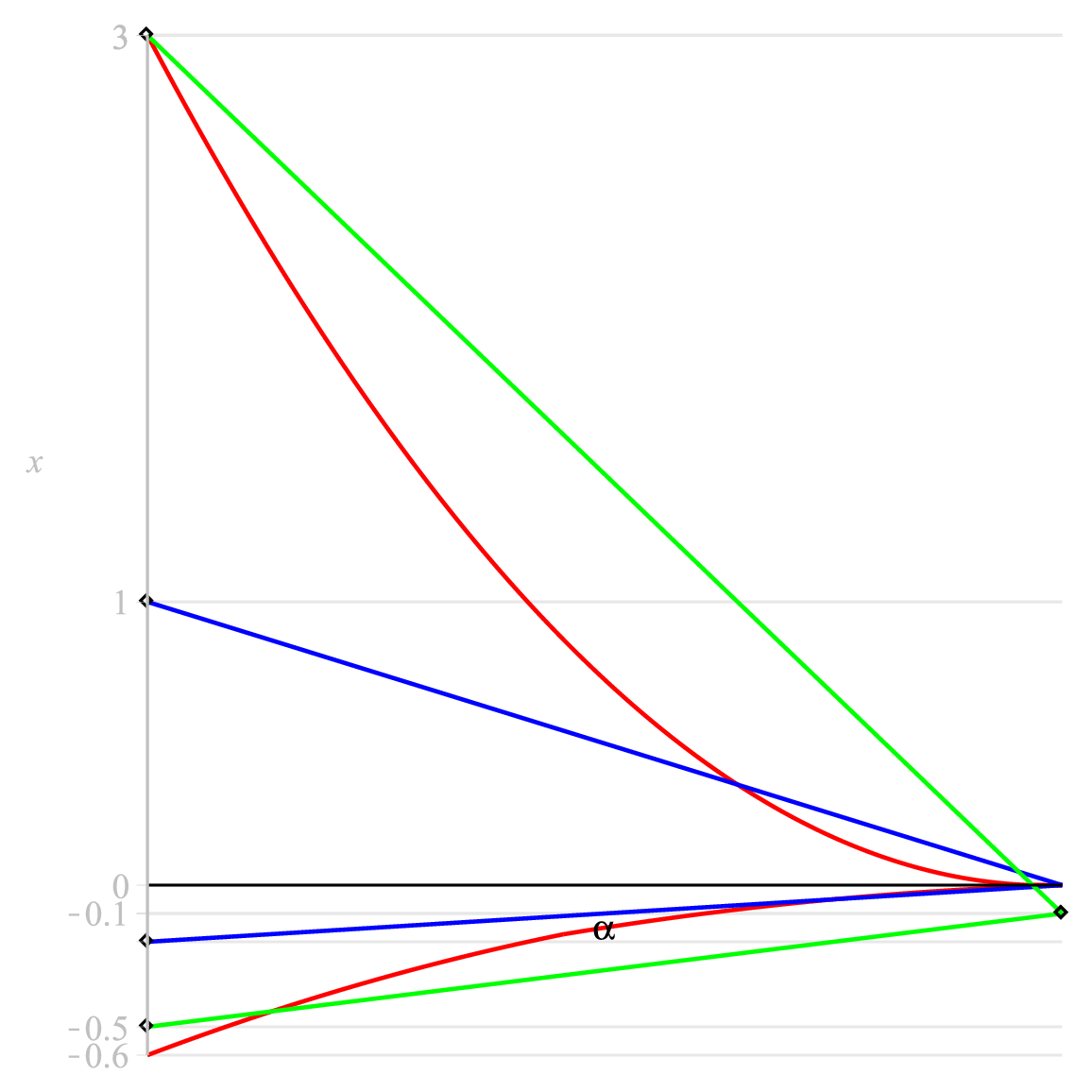}
\includegraphics[width=41mm]{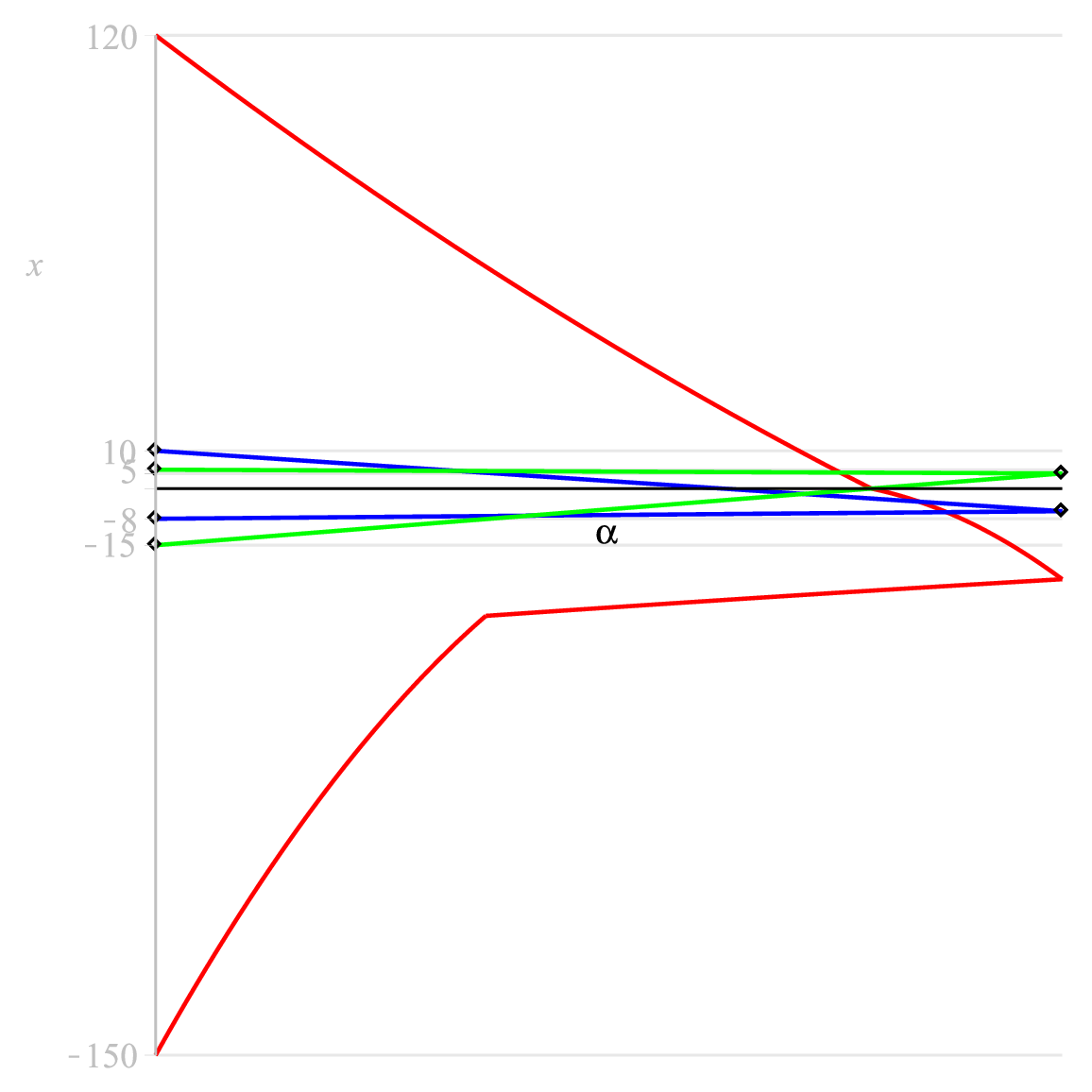}
\includegraphics[width=41mm]{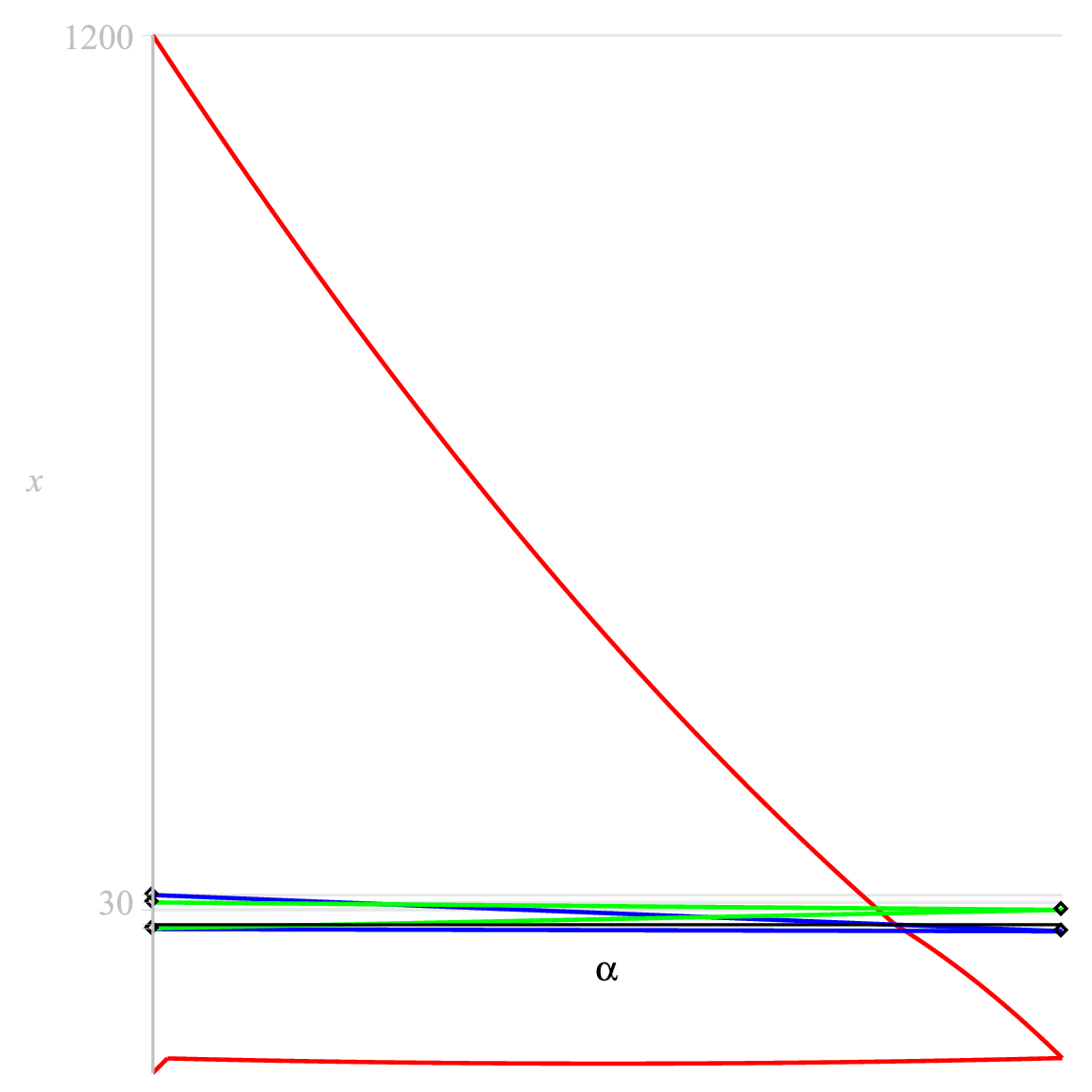}
\includegraphics[width=41mm]{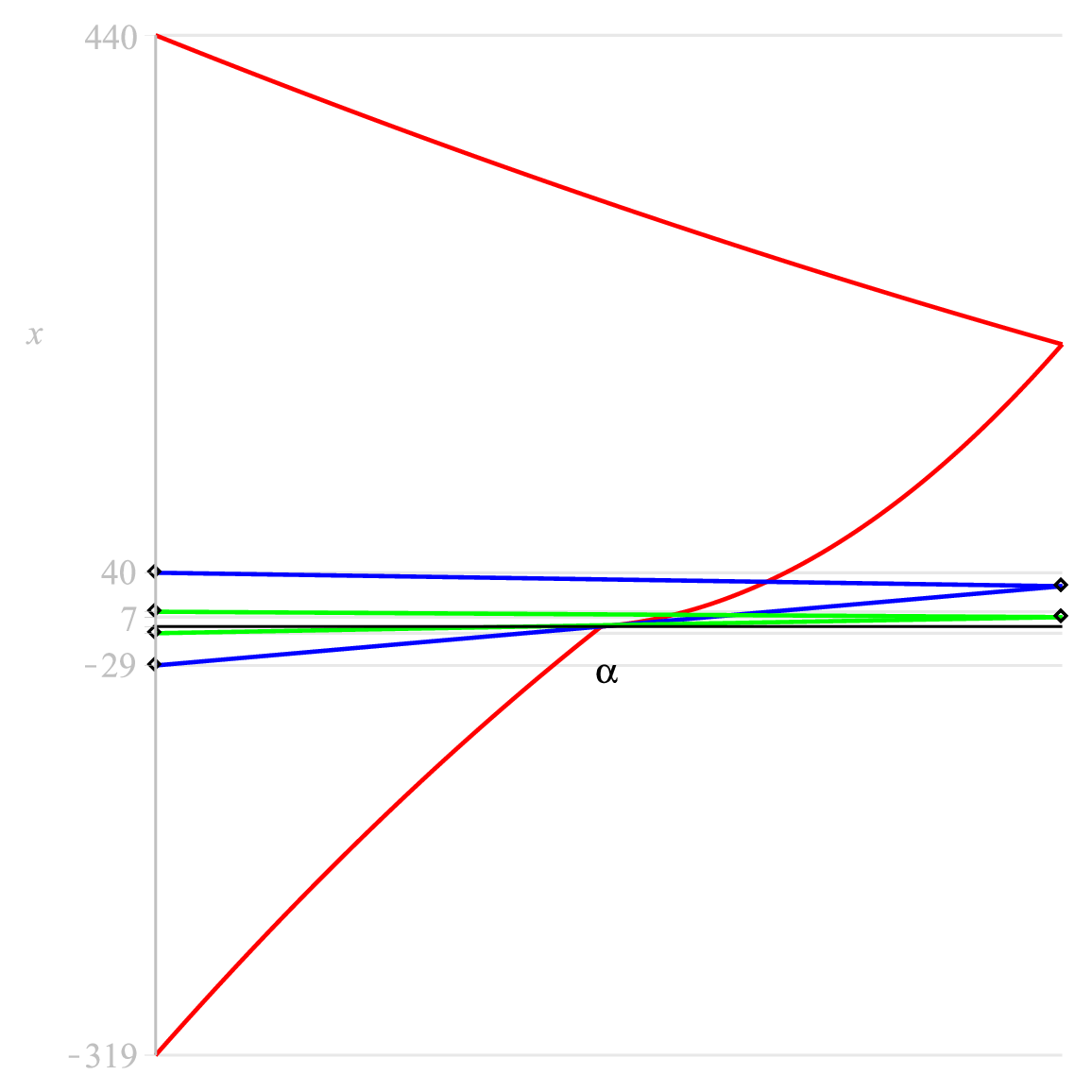}
\end{figure}

\vspace{1cm}
\begin{subsection}{Algebraic properties of the class $\Fth$ ({equivalently} $\Fth^{-1}$)}
\par
\smallskip As was the case for $\Fo$, $\Ft$ we have associativity and commutativity.\\
Compared to sections 1 and 2 we gain subtraction and division (Division could have been introduced from the outset in section 1 by~\ref{division3}). We realize however that neither really constitutes an inverse operation in the group sense to the basic operations of addition and multiplication.

\smallskip\noindent The distributive property of $\Fo$, $\Ft$ is lost. We only retain sub-distributivity in the sense that for every $\alpha\in[0,1]$ (that is for every $\alpha-$cut, at every $\alpha-$level) for three fuzzy intervals $[a_d(\alpha),a_u(\alpha)]$,
 $[b_d(\alpha),b_u(\alpha)]$ and $[c_d(\alpha),c_u(\alpha)]$,
\begin{equation}\label{distributivity3}
[a_d,a_u]\odot\bigl([b_d,b_u]\oplus[c_d,c_u]\bigr)\subseteq\bigl([a_d,a_u]\odot[b_d,b_u]\bigr)\oplus\bigl([a_d,a_u]\odot[c_d,c_u]\bigr)
\end{equation}
holds as sets.

\smallskip\noindent Graphically the above means that the graph of the left hand side of (\ref{distributivity3}) lies entirely within the graph of the right hand side.

\end{subsection}
\end{section}

\end{document}